\pdfoutput=1
\documentclass[final,a4paper,leqno,11pt]{article}
\usepackage{graphicx}
\usepackage{mathptmx}
\usepackage{amssymb,amsmath, amsfonts,amsthm}
\usepackage{a4wide}
\usepackage{color}
\def\sce{\setcounter{equation}{0}}
\newcommand{\red}[1]{\textcolor{red}{#1}}
\newcommand{\email}[1]{{\tt #1}}
\newcommand{\R}{\mathbb{R}}
\newcommand{\N}{\mathbb{N}}
\newcommand{\norm}[1]{\|#1\|}

\newcommand{\dist}[1]{{\rm dist}(#1)}
\newcommand{\B}{{\cal B}}
\newcommand{\K}{{\cal K}}
\newcommand{\M}{{\cal M}}
\newcommand{\J}{{\cal J}}
\newcommand{\F}{{\cal F}}
\newcommand{\Pp}{{\cal P}}

\newcommand{\mv}{\,\vert\,}

\newcommand{\longsetto}[1]{\mathop{\longrightarrow}\limits^#1}
\newcommand{\skalp}[1]{\langle #1\rangle}
\newcommand{\xb}{\bar x}
\newcommand{\yb}{\bar y}
\newcommand{\zb}{\bar z}
\newcommand{\pb}{\bar p}

\newcommand{\Lb}{\bar\Lambda}

\newcommand{\xba}{{\bar x^\ast}}

\newcommand{\oo}{o}
\newcommand{\OO}{{O}}

\newcommand{\argmax}{\mathop{\rm arg\,max}\limits}

\newcommand{\lin}{{\rm lin\,}}
\newcommand{\cl}{{\rm cl\,}}
\newcommand{\ri}{{\rm ri\,}}

\newcommand{\gph}{\mathrm{gph}\,}
\newcommand{\tto}{\rightrightarrows}
\def\h{\hfill\Box}
\def\disp{\displaystyle}
\def\dn{\downarrow}
\def\la{\langle}
\def\ra{\rangle}
\def\st{\stackrel}
\def\O{\Omega}
\def\ve{\varepsilon}
\def\Hat{\widehat}
\def\ox{\bar{x}}
\def\op{\bar{p}}
\def\oy{\bar{y}}
\def\oz{\bar{z}}
\def\ov{\bar{v}}
\def\ou{\bar{u}}

\def\gph{\mbox{\rm gph}\,}
\def\dist{\mbox{\rm dist}\,}

\def\ker{\mbox{\rm ker}\,}

\def\rge{{\rm rge\,}}

\def\lm{\lambda}

\def\kk{\kappa}
\def\Th{\Theta}
\def\Tilde{\widetilde}
\def\tilde{\widetilde}
\def\th{\theta}
\def\ph{\varphi}
\def\hat{\Hat}

\newlength{\myparboxwidth}
\newtheorem{Theorem}{Theorem}[section]
\newtheorem{proposition}[Theorem]{Proposition}

\newtheorem{lemma}[Theorem]{Lemma}
\newtheorem{corollary}[Theorem]{Corollary}

\newtheorem{example}[Theorem]{Example}

\title{Second-Order Variational Analysis\\of Parametric Constraint and Variational Systems}
\author{HELMUT GFRERER\footnote{Institute of Computational Mathematics, Johannes Kepler University Linz, A-4040 Linz, Austria; \email{helmut.gfrerer@jku.at}}
$\;$ \and $\;$ BORIS S. MORDUKHOVICH\footnote{Department of Mathematics, Wayne State University, Detroit, MI 48202, USA, and RUDN University, Moscow 117198, Russia; \email{boris@math.wayne.edu}}}
\begin{document}
\maketitle\vspace*{-0.3in}

\small
{\bf Abstract.} This paper is devoted to the generalized differential study of the normal cone mappings associated with a large class of {\em parametric constraint systems} (PCS) that appear, in particular, in nonpolyhedral conic programming. Conducting a local second-order analysis of such systems, we focus on computing the (primal-dual) {\em graphical derivative} of the {\em normal cone mapping} under the {\em $C^2$-cone reducibility} of the constraint set together with the fairly weak {\em metric subregularity constraint qualification} and its uniform parametric counterpart known as {\em Robinson stability}. The obtained precise formulas for computing the underlying {\em second-order} object are applied to the derivation of comprehensive conditions ensuring the important stability property of {\em isolated calmness} for solution maps to parametric variational systems associated with the given PCS.\vspace*{0.05in}

\vspace*{0.05in}

\vspace*{0.05in}

{\bf Abbreviated title.} Second-order variational analysis
\vspace*{-0.15in}

\normalsize
\section{Introduction and Initial Discussions}\sce\vspace*{-0.1in}

Our starting point here is a general class of {\em parametric constraint systems} (PCS) given by
\begin{equation}\label{EqConstrSystem}
g(p,x)\in C\;\mbox{ with }\;x\in X:=\R^n\;\mbox{ and }\;p\in P:=\R^m,
\end{equation}
where $x$ is the {\em decision} variable and $p$ is the {\em perturbation parameter}. Imposing the standard smoothness assumption on $g\colon P\times X\to Y$ with $Y:=\R^l$ and the $C^2$-{\em cone reducibility} assumption on $C\subset Y$ (see below), we define the set-valued {\em solution map} $\Gamma\colon P\rightrightarrows X$ to \eqref{EqConstrSystem} by
\begin{eqnarray}\label{EqGamma}
\Gamma(p):=\big\{x\in X\big|\;g(p,x)\in C\big\}\;\mbox{ for all }\;p\in P
\end{eqnarray}
and fix the reference feasible pair $(\pb,\xb)\in\gph\Gamma$ from the graph of $\Gamma$. Furthermore, we associate with \eqref{EqGamma} the {\em normal cone mapping} defined by
\begin{equation}\label{EqNormalConeMap}
\Psi(p,x):=N_{\Gamma(p)}(x)\;\mbox{ for all }\;p\in P\;\mbox{ and }\;x\in\Gamma(p),
\end{equation}
where the normal cone $N_{\Gamma(p)}$ is understood in the classical sense of convex analysis provided that the sets $\Gamma(p)$ are convex for all $p$. However, we do not impose any convexity of $\Gamma(p)$ in this paper, while our assumptions throughout the paper ensure that all the major normal cones of variational analysis (regular, limiting, and convexified) agree for the sets under consideration, and thus we may use the generic normal cone symbol $N$ below; see Section~2 for more details.

In what follows we fix a generalized normal $\ox^*\in N_{\Gamma(\pb)}(\xb)$ and pay the main attention to computing the {\em graphical derivative} $D\Psi(\ox,\op,\ox^*)$ of the normal cone mapping \eqref{EqNormalConeMap} at $(\ox,\op,\ox^*)$, which is a {\em primal-dual} object of second-order variational analysis defined and discussed in Section~2. Such a computation was first done by Rockafellar and Wets \cite{RoWe98} when $C$ is a convex polyhedron and \eqref{EqConstrSystem} satisfies the metric regularity condition, which is equivalent to the Mangasarian-Fromovitz qualification condition (MFCQ) in the polyhedral/nonlinear programming (NLP) setting under consideration. A number of recent publications has been devoted to computing $D\Psi$ in some particular cases of {\em nonparameterized} constraint systems \eqref{EqConstrSystem}, i.e., when $g=g(x)$ therein, under weaker assumptions on $g$ and $C$ with further applications to stability and optimization; see \cite{CH17,CHN17,GM15,GO16,GO17,HMS17,MOR15,MOR15a} for precise results and discussions. We specially emphasize the approach by Gfrerer and Outrata \cite{GO16} on computing $D\Psi$ in the {\em polyhedral} setting of $C$ without assuming MFCQ and/or other standard qualification conditions. The approach of \cite{GO16} was then employed by Chieu and Hien \cite{CH17} with imposing merely the {\em metric subregularity constraint qualification} (MSCQ) on the NLP constraint system. The quite recent paper by Hang, Mordukhovich and Sarabi \cite{HMS17} is the first one in the literature to compute $D\Psi$ by developing a novel approach that relays only on the MSCQ condition in the {\em nonpolyhedral} setting of \eqref{EqConstrSystem} with $g=g(x)$ generated by the second-order/Lorentz/ice-cream cone $C$ therein.

The main goal of this paper is to compute $D\Psi$ for the general parametric constraint systems \eqref{EqConstrSystem} with establishing significantly new results even for the {\em nonparametric setting} under imposing merely the MSCQ condition. Indeed, the precise formulas derived below are applied to PCS \eqref{EqConstrSystem} generated by {\em any $C^2$-cone reducible set} $C$, which covers the most interesting theoretically and important in applications constraint systems in conic programming including those generated by products of the ice-cream cones, the SDP cone of symmetric positive semidefinite matrices, etc. {\em No other assumptions but MSCQ} are required in the nonparametric case of \eqref{EqConstrSystem}.

The general {\em parametric setting} of \eqref{EqConstrSystem} is more involved and diverse. We impose here a {\em uniform} parametric counterpart of MSCQ called the {\em Robinson stability} (RS) in \cite{GfrMo17a}. If $g$ does not depend on the parameter, the Robinson stability reduces to MSCQ, whereas in case when
$\nabla_p g(\op,\ox)$ is surjective, RS is equivalent to the {\em metric regularity} qualification condition. Under RS we obtain lower and upper inclusions for the graphical derivative $D\Psi(\ox,\op,\ox^*)$. These estimates are shown to be exact under a certain additional assumption. If $g(p,x)$ {\em weakly} depends on the parameter ($\nabla_p g(\op,\ox)=0$), this assumption is automatically fulfilled, and we are in the same framework as in the nonparametric setting. In the case where $g$ {\em strongly} depends on the parameter (as, e.g., when $\nabla_p g(\op,\ox)$ is surjective), then this assumption amounts to a nondegeneracy condition.
\if{
 If $g(p,x)$ {\em weakly} depends on the parameter ($\nabla_p g(\op,\ox)=0$) at the reference point, the Robinson stability reduces to MSCQ, and we are in the same framework as in the nonparametric setting. In the case where $g$ {\em strongly} depends on the parameter (as, e.g., when $\nabla_p g(\op,\ox)$ has full rank), RS is close to the {\em metric regularity} qualification condition.}\fi

One alternative to the latter is the {\em polyhedrality} in the reduction, while we discover another one, which is the {\em strict complementarity} condition in the vein of Bonnans and Shapiro \cite{BonSh00}.

The obtained computation formulas for $D\Psi$ have high potentials for various applications to optimization, equilibrium, and related problems. We focus here on applications to {\em isolated calmness} of solutions maps for {\em parametric variational systems} (PVS) associated with the normal cone mapping \eqref{EqNormalConeMap}. Employing these formulas allows us to establish verifiable conditions for the validity of a new {\em H\"olderian} counterpart of isolated calmness for the PVS under consideration and then derive necessary and sufficient conditions for isolated calmness (Lipschitzian in nature) in both cases of weak and strong perturbations discussed above. It seems that the results obtained here are the best and most general in the literature even for parameter-independent mappings $\Gamma$ in \eqref{EqGamma}.

The rest of the paper is organized as follows. In Section~2 we recall and briefly discuss some notions from variational analysis and generalized differentiation that are systematically used in the sequel. Our standing assumptions are also formulated and discussed in this section together with preliminary results employed below. Section~3 studies second-order properties of PVS under the MSCQ assumption imposed on the reduction system. In Section~4 we derive formulas for computing the graphical derivative in the reduced PCS under imposing the RS condition on the reduction data.

The subsequent Section~5 contains major results on computing $D\Psi$ in terms of the given data of PCS \eqref{EqConstrSystem}. These results are applied in Section~6 to deriving verifiable conditions for isolated calmness of parametric variational systems associated with \eqref{EqConstrSystem}. We also present here a numerical example that illustrates the significance of the (rather mild) assumptions made for deriving the computation formulas and the obtained conditions for isolated calmness. Section~7 summarizes the main paper developments and discusses some directions of the future research.

Throughout the paper we use the standard notation of variational analysis (see, e.g., \cite{M06,RoWe98}) together with more special ones defined in the text. Recall that $\N:=\{1,2,\ldots\}$, $\B_X$ is the closed unit ball of $X$, $[x]:=\{tx\mv t\in\R\}$, and $\lin K$ is the largest linear subspace contained in a closed convex cone $K$; it is known as the {\em lineality space} of $K$. Observe that $\lin K=K\cap(-K)$, and that its annihilator $(\lin K)^\perp=K^\ast-K^\ast$ is the smallest subspace containing the dual cone $K^\ast$. To ease the reading, we want to clearly distinguish between primal and dual objects in our notation. Dual objects are marked by an asterisk except multipliers that are denoted by the Greek letters $\lambda$ and $\mu$ for historical reasons. Although all spaces $X,P,Y\ldots$ appearing in the paper are finite-dimensional, we do not identify the dual spaces $X^\ast,P^\ast,Y^\ast,\ldots$ with the primal ones. For differentiable mappings $f\colon U\to V$, the (Fr\'echet) derivative at a point $\bar u$ is denoted by $\nabla f(\bar u)$ and is a linear mapping from $U$ to $V$. Consequently, its adjoint mapping $\nabla f(\bar u)^\ast$ is a mapping between the dual spaces $V^\ast$ and $U^\ast$. Bidual spaces are identified with the primal spaces. E.g., given a mapping $f\colon P\times X\to X^\ast$, its derivative at a point $(\bar p,\bar x)$ is a linear mapping from $P\times X$ to $X^\ast$ and its adjoint $\nabla f(\bar p,\bar x)^\ast$ is a linear mapping from $X$ to $P^\ast\times X^\ast$.
\vspace*{-0.15in}

\section{Standing Assumptions and Preliminaries}\sce\vspace*{-0.1in}

We start with recalling the notion of $C^2$-cone reducibility of sets that has been well recognized and employed in optimization. Following \cite{BonSh00}, a closed set $C\subset Y$ is said to be {\em $C^2$-reducible} to the closed {\em cone} $K\subset E$ in finite dimensions at a point $\oy\in C$ if there exists a neighborhood $V\subset Y$ of $\yb$ and a twice continuously differentiable mapping $h\colon V\to E$ \red{}
such that
\begin{eqnarray}\label{DefConeReducible}
h(\yb)=0,\;\nabla h(\yb)\;\mbox{ is surjective },\;\mbox{ and }\;C\cap V=\big\{y\in V\big|\;h(y)\in K\big\}.
\end{eqnarray}
As discussed in \cite{BonSh00}, the most important constraint sets in conic programming (including products of the ice-cream cones, the SDP cone, etc.) are $C^2$-cone reducible at all their feasible points $\oy\in C$.

The next basic notion systematically used in the paper is a weak qualification condition coined in \cite{GM15} as the {\em metric subregularity constraint qualification} (MSCQ). We say that MSCQ holds for the system $g(z)\in C$ at the feasible point $\oz$ if the set-valued mapping $M(z):=g(z)-C$ is {\em metrically subregular} at $(\zb,0)$, i.e., there exists a constant/modulus $\kk\ge 0$ such that
\begin{eqnarray}\label{DefMetrSubregCQ}
{\rm dist}\big(z;M^{-1}(0)\big)\le\kk\,{\rm dist}\big(0;M(z)\big)\;\mbox{ for all }\;z\;\mbox{ close to }\;\oz.
\end{eqnarray}
The reader is referred to \cite{GM15} for more details on this constraint qualification, sufficient conditions for its validity, and its relationships with other qualification conditions in NLPs. Note that MSCQ is essentially weaker than its {\em metric regularity} counterpart, which corresponds to \eqref{DefMetrSubregCQ} with replacing $\oy=0$ in both parts therein by any $y\in Y$ of small norm. The latter condition (but not MSCQ) reduces to MFCQ and the Robinson constraint qualification in NLPs and conic programming, respectively.

A crucial role in the general parametric setting of \eqref{DefMetrSubregCQ} is played by the stability notion introduced recently in \cite{GfrMo17a} as an extension of the one first studied by Robinson \cite{Rob76} in particular situations. We say that PCS \eqref{EqConstrSystem} enjoys the {\em Robinson stability} (RS) property at $(\pb,\xb)$ with modulus $\kappa\ge 0$ if there are neighborhoods $U$ of $\xb$ and $W$ of $\pb$ such that
\begin{equation}\label{EqRS}
{\rm dist}\big(x;\Gamma(p)\big)\le\kappa\,{\rm dist}\big(g(p,x);C\big)\;\mbox{ for all }\;(p,x)\in W\times U
\end{equation}
in terms of the usual point-to-set distance. This property was largely studied in \cite{GfrMo17a} with deriving efficient first-order and second-order conditions for its validity and various applications. It is clear from \eqref{EqRS} that RS can be interpreted as metric subregularity of the mapping $x\mapsto g(p,x)-C$ at $(\ox,0)$ for every point $x\in\Gamma(p)$ near $\ox$, where the modulus $\kk$ is {\em uniform} with respect to parameters $p\in W$.

For the rest of the paper we impose the following requirements on the given data of \eqref{EqConstrSystem}.\\[1ex]
{\bf Standing Assumptions:} {\bf (i)} The mapping $g\colon P\times X\to Y$ is $C^2$-smooth around $(\op,\ox)$.\\
{\bf (ii)} The set $C\subset Y$ is $C^2$-reducible at $g(\zb)$ to a closed and convex cone $K\subset E$.\vspace*{0.05in}

Besides these standing assumptions, each result below requires the validity of either the MSCQ or RS condition at the reference point.\vspace*{0.05in}

Now it is time for entering and discussing constructions of generalized differentiation widely used in the paper. First we clarify the meaning of the normal cone in \eqref{EqGamma} in the nonconvex setting. Given a closed set $\O\subset Z$ with $\oz\in\O$, the (Fr\'echet) {\em regular normal cone} to $\O$ at $\oz$ is defined by
\begin{eqnarray}\label{rn}
\Hat N_{\O}(\oz):=\disp\Big\{z^*\in Z^*\Big|\;\limsup_{z\st{\O}{\to}\oz}\frac{\la z^*,z-\oz\ra}{\|z-\oz\|}\le 0\Big\},
\end{eqnarray}
where the symbol $z\st{\O}{\to}\oz$ indicates that $z\to\oz$ with $z\in\O$. It has been realized that \eqref{rn} fails to have good properties required for generalized normals to nonconvex sets, namely: it may be trivial $(=0)$ at boundary points, nonrobust with respect to perturbations of the initial data, does not enjoy calculus rules, etc. As shown in \cite{M06,RoWe98}, the aforementioned and other desired properties hold--in spite of its nonconvexity--for the (Mordukhovich) {\em limiting normal cone} to $\O$ at $\oz\in\O$ defined by
\begin{eqnarray}\label{ln}
N_\O(\oz):=\big\{z^*\in Z^\ast\big|\;\exists\,z_k\st{\O}{\to}\oz,\;z^*_k\in\Hat N_\O(z_k)\;\mbox{ with }\;z^*_k\to z^*\;\mbox{ as }\;k\to\infty\big\},
\end{eqnarray}
which is however more challenging in computation. Luckily, the assumptions imposed in this paper ensure that the normal cones \eqref{rn} and \eqref{ln} agree for the sets under considerations, and thus we can combine good properties of $N_\O$ with the computational advantages of $\Hat N_\O$. Indeed, the $C^2$-reducibility \eqref{DefConeReducible} of the set $C$ at $\oy\in C$ to the convex cone $K$ yields by the surjectivity of the derivative $\nabla h(\oy)$ and the well-known calculus rule for normals to inverse images (see, e.g., \cite[Theorem~6.14]{RoWe98}) that
\begin{eqnarray}\label{nr1}
\Hat N_C(\oy)=\Hat N_{h^{-1}(K)}\big(h(\oy)\big)=\nabla h(\oy)^*\Hat N_K\big(h(\oy)\big)=\nabla h(\oy)^*N_K\big(h(\oy)\big)=N_C(\oy).
\end{eqnarray}
Furthermore, the MSCQ (and hence RS) conditions imposed on \eqref{EqConstrSystem} at $(\op,\ox)$ imply that
\begin{eqnarray}\label{nr2}
\Hat N_{g^{-1}(C)}(\op,\ox)=\nabla g(\op,\ox)^*\Hat N_C\big(g(\op,\ox)\big)=\nabla g(\op,\ox)^*N_C\big(g(\op,\ox)\big)=N_{g^{-1}(C)}(\op,\ox),
\end{eqnarray}
which follows from \cite[Theorem~4.1]{HJO02} and \eqref{nr1}. Thus we have by \eqref{EqGamma} that
\begin{eqnarray}\label{nr3}
\Hat N_{{\rm gph}\,\Gamma}(\op,\ox)=N_{{\rm gph}\,\Gamma}(\op,\ox)\;\mbox{ and }\;\Hat N_{\Gamma(\op)}(\ox)=N_{\Gamma(\op)}(\ox).
\end{eqnarray}
The {\em normal regularity} relationships in \eqref{nr1}--\eqref{nr3} show that the normal cones therein agree also with the (Clarke) {\em convexified normal cone} to these sets, which is the convex closure of the limiting one \eqref{ln}. It allows us to use in this paper the generic symbol $N$ for all the three major cones.

Following \cite{BonSh00}, we say that the {\em strict complementarity} condition holds for the system $g(z)\in C$ at $(\oz,z^*)$ with $g(\oz)\in C$ and $z^*\in N_{g^{-1}}(\oz)$ if there is $\mu\in{\rm ri}N_{C}(g(\oz))$ such that $z^*=\nabla g(\oz)^*\mu$. We refer the reader to \cite{BonSh00} for more discussions on this condition and its relationships with the classical notions.

In conventional understanding, normals are dual-space constructions while tangents are primal-space ones. The tangent collection intensively used below is known as the (Bouligand-Severi) {\em tangent} or {\em contingent cone} to $\O$ at $\oz\in\O$ being defined by
\begin{eqnarray}\label{tan}
T_\O(\oz):=\big\{w\in Z\big|\;\exists\,t_k\dn 0,\;w_k\to w\;\mbox{ with }\;\oz+t_k w_k\in\O\;\mbox{ for all }\;k\in\N\big\}.
\end{eqnarray}
A tangent vector $w\in T_\O(\oz)$ is {\em derivable} if there exists a vector function $\xi\colon[0,\ve]\to\O$ with $\ve>0$ such that $\xi(0)=\oz$ and $\xi'_+(0)=w$, where $\xi'_+(\cdot)$ stands for the right derivative of $\xi(\cdot)$. By the finite dimensionality of the spaces in question we always have $\widehat N_\O(\oz)=T_\O(\oz)^*$. If the set $\O$ is normally regular at $\oz$, then there is the following {\em tangent-normal duality}:
\begin{eqnarray}\label{tan1}
T_\O(\oz)=N_\O(\oz)^*:=\big\{w\in Z\big|\;\la w,z^*\ra\le 0\;\mbox{if}\;z^*\in N_\O(\oz)\big\}\;\mbox{ and }\;N_\O(\oz)=T_\O(\oz)^*.
\end{eqnarray}
Using the tangent cone \eqref{tan}, define next the {\em critical cone} to $\O$ at $(\oz,\oz^*)$ with $\oz\in\O$ and $\oz^*\in N_\O(\oz)$, which is generated by \eqref{tan} and the orthogonal complement $\{\oz^*\}^\perp:=\{w\in Z|\;\la w,\oz^*\ra=0\big\}$ of $\oz^*$ as
\begin{eqnarray}\label{cc}
\K_\Omega(\zb,\oz^*):=T_\Omega(\zb)\cap\big\{\oz^*\}^\perp.
\end{eqnarray}
Another primal-space construction generated by the tangent cone \eqref{tan} is the {\em graphical derivative} of a set-valued mapping $\Psi\colon Z\tto U$ at the given point $(\oz,\ou)\in\gph\Psi$ defined by
\begin{eqnarray}\label{gd}
D\Psi(\oz,\ou)(w):=\big\{v\in U\big|\;(w,v)\in T_{{\rm gph}\,\Psi}(\oz,\ou)\}\;\mbox{ for all }\;w\in Z
\end{eqnarray}
via the tangent cone \eqref{tan} to the graph of $\Psi$ at $(\oz,\ou)$. Our major goal in what follows is to compute the graphical derivative \eqref{gd} for the normal cone mapping $\Psi$ given in \eqref{EqNormalConeMap}.

We conclude this section with the following auxiliary proposition, which combines some facts often used below in the proofs of the main results of the paper.\vspace*{-0.07in}

\begin{proposition}[\bf normal and critical directions to inverse images]\label{LemNormalCone} Let $f\colon Z\to Y$ be $C^1$-smooth around $\oz$ with $f(\oz)\in\O$, where $\O$ is a closed convex set. If the constraint system $f(z)\in\O$ satisfies MSCQ at $\oz$ with modulus $\kk$, then
$N_{f^{-1}(\O)}(\zb)=\nabla f(\zb)^\ast N_\O(g(\zb))$ and for every $z^*\in N_{f^{-1}(\O)}(\zb)$ we have
\begin{equation}\label{EqBndMult}
z^*\in\nabla f(\zb)^\ast\big(N_\O(f(\zb))\cap\kappa\norm{z^*}\B_{Y^*}\big).
\end{equation}
Furthermore, for all such $z^*\in N_{f^{-1}(\O)}(\zb)$ the dual cone to \eqref{cc} is represented by
\begin{eqnarray}\label{cc1}
\K_{f^{-1}(\O)}(\zb,z^*)^\ast=\cl\big(\nabla f(\zb)^\ast N_\O(f(\zb))+[z^*]\big),
\end{eqnarray}
where the closure operation can be omitted if the strict complementarity condition holds at $(\zb,z^*)$.
\end{proposition}\vspace*{-0.1in}
{\bf Proof.} Inclusion \eqref{EqBndMult} follows from \cite[Lemma~2.1]{GfrMo17a} and implies that $N_{f^{-1}(\O)}(\zb)\subset\nabla f(\zb)^\ast N_\O(f(\zb)$. The opposite inclusion is a consequence of \cite[Theorem~6.14]{RoWe98}, and thus we get the claimed equality for $N_{f^{-1}(\O)}(\zb)$. To proceed with \eqref{cc1}, pick any
$z^*\in N_{f^{-1}(\O)}(\zb)$ and observe that
\begin{eqnarray*}
\K_{f^{-1}(\O)}(\zb,z^*)^\ast=\big(T_{f^{-1}(\O)}(\zb)\cap\{z^*\}^\perp\big)^\ast=\cl\big(N_{f^{-1}(\O)}(\zb)+[z^*]\big)=\cl\big(\nabla f(\zb)^\ast N_\O(f(\zb))+[z^*]\big),
\end{eqnarray*}
which verifies \eqref{cc1}. If strict complementarity holds at $(\zb,z^*)$, take $\mu\in\ri N_\O(f(\zb))$ with $z^*=\nabla f(\zb)^\ast\mu$ and get by the convexity of $N_\O(f(\zb))$ that $N_\O(f(\zb))+\R_+\mu\subset N_\O(f(\zb))$. This yields
\begin{eqnarray*}
\nabla f(\zb)^\ast N_\O\big(f(\zb)\big)+[z^*]&=&\nabla f(\zb)^\ast\big(N_\O(f(\zb))+[\mu]\big)=
\nabla f(\zb)^\ast\big(N_\O(f(\zb))-\R_+\mu\big)\\
&=&\nabla f(\zb)^\ast\big(\R_+\big(N_\O(f(\zb))-\mu)\big)=\cl\Big(\nabla f(\zb)^\ast\big(\R_+(N_\O(f(\zb))-\mu\big)\Big)\\
&=&\cl\big(\nabla f(\zb)^\ast N_\O(f(\zb))+[z^*]\big)
\end{eqnarray*}
since $\R_+(N_\O(f(\zb))-\mu)$ is the subspace parallel to the affine hull of $N_\O(f(\zb))$ and hence the set $\nabla f(\zb)^\ast(\R_+(N_\O(f(\zb))-\mu))$ is a subspace as well. The closedness of the later shows that the closure operation can be removed from \eqref{cc1}, which completes the proof of the proposition. $\h$\vspace*{-0.15in}

\section{Second-Order Properties of PCS under Reduction}\sce\vspace*{-0.1in}

In this section we study, under the standing assumptions formulated above, some second-order properties on the {\em reduced system} with imposing in addition the {\em MSCQ} condition on its data. Denoting $\oz:=(\op,\ox)\in Z:=P\times X$ and using the reduction data $h$ and $K$ from \eqref{DefConeReducible}, we consider the mapping
\begin{equation}\label{EqMapG}
G\colon Z\to E\;\mbox{ with }\;G(z):=h\big(g(z)\big)\;\mbox{ and }\;G(\zb)=0,
\end{equation}
where the latter condition is essential in what follows. Observe that the reduced PCS given by $G(p,x)\in K$ is {\em equivalent} to the original one \eqref{EqConstrSystem} in the following sense. Taking the neighborhood $V$ of $g(\pb,\xb)$ from \eqref{DefConeReducible} and choosing open neighborhoods $U$ of $\xb$ and $Q$ of $\pb$ with $Q\times U\subset g^{-1}(V)$, for each $p\in Q$ we get that $\Gamma(p)\cap U=\Tilde\Gamma(p)\cap U$, where $\Tilde\Gamma(p):=\{x\mv G(p,x)\in K\}$. Therefore we have
$$
\widehat N_{\Gamma(p)}(x)=\widehat N_{\Gamma(p)\cap U}(x)=\widehat N_{\Tilde\Gamma(p)\cap U}(x)=\widehat N_{\Tilde\Gamma(p)}(x)
$$
whenever $(p,x)\in Q\times U$. This shows that the graphical derivative of the normal cone mapping $\Psi$ at $(\pb,\xb,\ov)$ coincides with the one of the mapping $(p,x)\rightrightarrows\widehat N_{\Tilde\Gamma(p)}(x)$. To simplify the notation, suppose that $\Gamma(p)=\Tilde\Gamma(p)$ for all $p$ sufficiently close to $\pb$. Note that $\gph\Gamma=\{z\in Z\mv g(z)\in C\}$, and hence $\gph\Gamma\cap W=\{z\in W\mv G(z)\in K\}$ for some neighborhood $W$ of $\zb$.

In the rest of this section we present three statements concerning second-order properties of the reduced system, which are of their own interest while playing an important role in establishing the main results of the paper. The first proposition justifies the (Robinson) {\em upper Lipschitzian} property (see, e.g., \cite{BonSh00,RoWe98}) of an auxiliary set-valued mapping built upon the data in \eqref{DefConeReducible} and \eqref{EqMapG}.\vspace*{-0.07in}

\begin{proposition}[\bf upper Lipschitzian property of the second-order auxiliary mapping]\label{LemUpperLip} Assume that MSCQ holds for the system $G(z)\in K\subset E$ at $\zb$ and take $v\in Z$ with $\nabla G(\zb)v\in K$. Then we can find a positive number $\kappa$ such that for every pair $(w,q)\in Z\times E$ with
\begin{equation}\label{EqUpperLip_wq}
\nabla G(\zb)w+\frac 12\nabla^2G(\zb)(v,v)+q\in T_K\big(\nabla G(\zb)v\big)
\end{equation}
there exists a vector $\Tilde w\in Z$ satisfying the conditions
\begin{equation}\label{EqUpperLip}
\nabla G(\zb)\tilde w+\frac 12\nabla^2G(\zb)(v,v)\in T_K(\nabla G(\zb)v)\;\mbox{ and }\;\norm{\tilde w-w}\le\kappa\norm{q}.
\end{equation}
The latter condition can be reformulated as the upper Lipschitzian property
\begin{eqnarray}\label{RemUpperLip}
\Theta(q)\subset\Theta(0)+\kappa\norm{q}\B_Z\;\mbox{ for all }\;q\in E
\end{eqnarray}
of the auxiliary set-valued mapping $\Th\colon E\tto Z$ defined by
\begin{eqnarray*}
\Theta(q):=\Big\{w\in Z\Big|\;\nabla G(\zb)w+\frac 12\nabla^2G(\zb)(v,v)+q\in T_K(\nabla G(\zb)v)\Big\}.
\end{eqnarray*}
\end{proposition}\vspace*{-0.05in}
{\bf Proof.} The imposed MSCQ gives us a neighborhood $W$ of $\zb$ and a number $\kappa>0$ such that ${\rm dist}(z;G^{-1}(K))\le\kappa{\rm dist}(G(z);K)$ for all $z\in W$. Consider $(w,q)$ satisfying \eqref{EqUpperLip_wq}, we observe that
\begin{eqnarray*}
T_K\big(\nabla G(\zb)v\big)=\cl\big(K+[\nabla G(\zb)v]\big)
\end{eqnarray*}
and fix $\epsilon>0$. Then there are $y_\epsilon\in K$, $q_\epsilon\in E$, and $\alpha_\epsilon\in\R$ with $\norm{q_\epsilon-q}\le\epsilon$ and $\nabla G(\zb)w+\frac 12\nabla^2G(\zb)(v,v)+q_\epsilon=y_\epsilon-\alpha_\epsilon\nabla G(\zb)(v)$. Given any positive number $t$, we have the representations
\begin{eqnarray*}
G\big(\zb+tv+t^2(w+\alpha_\epsilon v)\big)&=&t\nabla G(\zb)v+t^2\Big(\nabla G(\zb)(w+\alpha_\epsilon v)+\frac 12\nabla^2G(\zb)(v,v)\Big)+o(t^2)\\
&=&t\nabla G(\zb)v+t^2(y_\epsilon-q_\epsilon)+o(t^2).
\end{eqnarray*}
Since both $\nabla G(\zb)v$ and $y_\epsilon$ are contained in the convex cone $K$, we have $t\nabla G(\zb)v+t^2y_\epsilon\in K$, which yields the estimate $\dist(G(\zb+tv+t^2(w+\alpha_\epsilon v));K)\le t^2\norm{q_\epsilon}+o(t^2)$. Hence for every small number $t>0$ there is $s_\epsilon(t)$ with $\norm{s_\epsilon(t)}\le\kappa(\norm{q_\epsilon}+\frac{o(t^2)}{t^2})$ and $G\big(\zb+tv+t^2(w+\alpha_\epsilon v+s_\epsilon(t))\big)\in K$. Thus we can find a sequence $t_k\downarrow 0$ such that the sequence $s_\epsilon(t_k)$ converges to some $s_\epsilon$ satisfying $\norm{s_\epsilon}\le\kappa\norm{q_\epsilon}$. Using again the conic structure of $K$ tells us that the quantities
\begin{eqnarray*}
\frac{G\big(\zb+t_kv+t_k^2(w+\alpha_\epsilon v+s_\epsilon\big(t_k)\big)}{t_k^2}=\Big(\frac 1{t_k}+\alpha_\epsilon\Big)\nabla G(\zb)v+\nabla G(\zb)(w+s_\epsilon\big(t_k)\big)+\frac 12\nabla^2G(\zb)(v,v)+\frac{o(t_k^2)}{t_k^2}
\end{eqnarray*}
belong to the convex cone $K$. This clearly implies that
\begin{eqnarray*}
\nabla G(\zb)\big(w+s_\epsilon(t_k)\big)+\frac 12\nabla^2G(\zb)(v,v)+\frac{o(t_k^2)}{t_k^2}\in K+[\nabla G(\zb)v],
\end{eqnarray*}
which in turn yields by passing to the limit as $t\dn 0$ the inclusion
\begin{eqnarray*}
\nabla G(\zb)(w+s_\epsilon)+\frac 12\nabla^2G(\zb)(v,v)\in\cl\big(K+[\nabla G(\zb)v]\big)=T_K\big(\nabla G(\zb)v\big).
\end{eqnarray*}
Considering an arbitrary sequence $\epsilon_k\downarrow 0$ and passing to a subsequence if necessary, suppose that $s_{\epsilon_k}$ converges to some $s\in Z$ and then obtain $\norm{s}\le\kappa\norm{q}$ and $\nabla G(\zb)(w+s)+\frac 12\nabla^2G(\zb)(v,v)\in T_K(\nabla G(\zb)v)$, which verifies \eqref{EqUpperLip} with $\tilde w=w+s$ and thus completes the proof of the proposition.$\h$\vspace*{0.03in}

Given  $z^*\in N_{{\rm gph}\,\Gamma}(\oz)$ and a critical direction $v\in\K_{{\rm gph}\,\Gamma}(\zb,z^*)$, we define now the following (dual) problem of {\em conic linear programming}:
\begin{eqnarray}\label{dual}
(D_{v,z^*})\qquad\sup_\mu\frac 12\big\la\mu,\nabla^2G(\zb)(v,v)\big\ra\;\mbox{ subject to }\;\mu\in K^\ast\;\mbox{ and }\;\nabla G(\zb)^\ast\mu=z^*.
\end{eqnarray}
The next proposition establishes an important property in conic programming called the ``approximate duality" in \cite{HMS17}, where it was obtained for second-order cone programs by using their specific features. Here we proceed in the general framework by exploiting MSCQ.\vspace*{-0.1in}

\begin{proposition}[\bf approximate duality in conic linear programming]\label{PropConicLinProg} Assume that MSCQ holds for the system $G(z)\in K$ at $\zb$ and consider the conic linear program \eqref{dual} generated by some $z^*\in N_{{\rm\small gph}\,\Gamma}(\oz)$ and $v\in\K_{{\rm gph}\,\Gamma}(\zb,z^*)$.
Then the set of optimal solution to $(D_{v,z^*})$ is nonempty and for every optimal solution $\hat\mu$ and for every $\epsilon>0$ there is some $w_\epsilon$ satisfying the conditions
\begin{eqnarray}\label{EqApprDual1}
{\rm dist}\Big(\nabla G(\zb)w_\epsilon+\frac 12\nabla^2G(\zb)(v,v);K\Big)\le\epsilon\;\mbox{ and }\;\skalp{z^*,w_\epsilon}+\frac 12\big\la\hat\mu,\nabla^2G(\zb)(v,v)\big\ra\ge-\epsilon.
\end{eqnarray}
\end{proposition}\vspace*{-0.1in}
{\bf Proof.} Consider the primal conic linear program $(\tilde P_{v,z^*})$ and its dual $(\tilde D_{v,z^*})$ given by
\begin{eqnarray*}
(\tilde P_{v,z^*})\qquad\inf_w-\skalp{z^*,w}\;\mbox{ subject to }\;\nabla G(\zb)w+\frac 12\nabla^2G(\zb)(v,v)\in T_K\big(\nabla G(\zb)v\big)\;\mbox{ and}
\end{eqnarray*}
\begin{eqnarray*}
(\tilde D_{v,z^*})\qquad\sup_\mu\frac 12\big\la\mu,\nabla^2G(\zb)(v,v)\big\ra\;\mbox{ subject to }\;\mu\in N_K\big(\nabla G(\zb)v\big),\;\nabla G(\zb)^\ast\mu=z^*,
\end{eqnarray*}
respectively. By the well-known relationship $N_K(\nabla G(\zb)v)=K^\ast\cap\big[\nabla G(\zb)v\big]^\perp$, for every $\mu$ feasible to $(D_{v,z^*})$ we get $\skalp{\mu,\nabla G(\zb)v}=\skalp{z^*,v}=0$, which shows that the programs $(\tilde D_{v,z^*})$ and $(D_{v,z^*})$ from \eqref{dual} are equivalent. It follows from the upper Lipschitzian property \eqref{RemUpperLip} verified in Proposition~\ref{LemUpperLip} and the results of \cite[Propositions~2.147 and 2.186]{BonSh00} that the set of optimal solutions to $(D_{v,z^*})$ is nonempty and that the optimal values of $(\tilde P_{v,z^*})$ and $(D_{v,z^*})$ agree. Hence for every optimal solution $\hat\mu$ to \eqref{dual} and every $\epsilon>0$ there is some $\hat w_\epsilon$ feasible to $(\tilde P_{v,z^\ast})$ satisfying $-\skalp{z^*,\hat w_\epsilon}-\frac 12\skalp{\hat\mu,\nabla^2G(\zb)(v,v)}\le\epsilon$. Furthermore, we can find a real number $\alpha_\epsilon$ such that
\begin{eqnarray*}
{\rm dist}\Big(\nabla G(\zb)\hat w_\epsilon+\frac 12\nabla^2G(\zb)(v,v)+\alpha_\epsilon\nabla G(\zb)v;K\Big)\le\epsilon.
\end{eqnarray*}
Taking finally into account that $\skalp{z^*,v}=0$ as shown above, we arrive at both conditions in \eqref{EqApprDual1} with $w_\epsilon=\hat w_\epsilon+\alpha_\epsilon v$ and thus complete the proof of the proposition. $\h$\vspace*{0.05in}

To proceed further, fix $z^*\in N_{{\rm gph}\,\Gamma}(\zb)$, form the set of multipliers
\begin{eqnarray*}
\M(\zb,z^*):=\big\{\mu\in K^\ast\big|\;\nabla G(\zb)^\ast\mu=z^*\},
\end{eqnarray*}
and define the multiplier set in the direction $v\in\K_{{\rm gph}\,\Gamma}(\zb,z^*)$ by
\begin{eqnarray}\label{mult}
\M(\zb,z^*;v):=\argmax\big\{\skalp{\mu,\nabla^2G(\zb)(v,v)}\mv\mu\in\M(\zb,z^*)\big\}.
\end{eqnarray}
Observe that $\M(\zb,z^*;v)$ is the set of optimal solutions to the dual program $(D_{v,z^*})$ from \eqref{dual}, and thus it is nonempty by Proposition~\ref{PropConicLinProg} under the imposed MSCQ condition.

Now we are ready to derive the main result of this section that establishes primal-dual relationships between second-order elements of the constraint system \eqref{EqConstrSystem} and of its reduction \eqref{EqMapG} under MSCQ. It is convenient to use in what follow the {\em product projection operator} defined by
\begin{eqnarray}\label{proj}
\pi_{X^*}\colon P^*\times X^*\to X^*\;\mbox{ with }\;\pi_{X^*}(p^*,x^*):=x^*,
\end{eqnarray}
where the decision space $X$ and the parameter space $P$ are taken from \eqref{EqConstrSystem}.\vspace*{-0.07in}

\begin{Theorem}[\bf second-order primal-dual relationships for PCS]\label{PropSuffTanDir} Assume that MSCQ holds for the system $G(z)\in K$ at $\zb$. Then for every $z^\ast\in N_{{\rm gph}\,\Gamma}(\zb)\cap\pi_{X^\ast}^{-1}(\xba)$, every $v=(q,u)\in\K_{{\rm gph}\,\Gamma}(\zb,z^\ast)$, every $\hat\mu\in\M(\zb,z^\ast;v)$, and every $w_{x^\ast}^\ast\in X^\ast$ such that there is a sequence of $w_k^\ast\in\K_{{\rm gph}\,\Gamma}(\zb,z^\ast)^\ast$ with
\begin{eqnarray*}
w_{x^\ast}^\ast=\lim_{k\to\infty}\pi_{X^\ast}(w_k^\ast)\;\mbox{ and }\;\lim_{k\to\infty}\skalp{w_k^\ast,v}=0
\end{eqnarray*}
we have the following inclusion into the tangent cone to the graph of the solution map $\Gamma$ from \eqref{EqGamma}:
\begin{eqnarray*}
\left(q,u,\nabla\big(\nabla_x G(\cdot)^\ast\hat\mu\big)(\zb)v+w_{x^\ast}^\ast\right)\in T_{{\rm gph}\,\Psi}(\pb,\xb,\xba),
\end{eqnarray*}
where these tangents are derivable. It gives us, in particular, that
\begin{eqnarray}\label{RemConv_w_k}
\left(q,u,\nabla\big(\nabla_x G(\cdot)^\ast\hat\mu\big)(\zb)v+\pi_{X^\ast}(w^\ast)\right)\in T_{{\rm gph}\,\Psi}(\pb,\xb,\xba)\;\mbox{ for all }\;w^\ast\in N_{\K_{{\rm gph}\,\Gamma}(\zb,z^\ast)}(v)
\end{eqnarray}
along with the elements $(z^*,v,u,q,\hat\mu)$ listed above.
\end{Theorem}\vspace*{-0.05in}
{\bf Proof.} Since $z^\ast=\nabla G(\zb)^\ast\hat\mu$ by \eqref{mult}, it follows from Proposition~\ref{LemNormalCone} that
\begin{eqnarray*}
\K_{{\rm gph}\,\Gamma}(\zb,z^\ast)^\ast=\cl\big(\nabla G(\zb)^\ast K^\ast+[z^\ast]\big)=\cl\big(\nabla G(\zb)^\ast(K^\ast+[\hat\mu])\big).
\end{eqnarray*}
Hence for every $k\in\N$ there exist $\tilde\mu_k\in K^\ast$ and $\alpha_k\in\R$ such that
$\norm{w_k^\ast-\nabla G(\zb)^\ast(\tilde\mu_k-\alpha_k\hat\mu)}\le k^{-1}$. Let us recursively construct a sequence of real numbers $\beta_k$ for $k\in\N$ by
\begin{eqnarray*}
\beta_1:=\max\big\{-\alpha_1+1,0\big\},\;\beta_{k+1}:=\max\big\{(k+1)(\vert\alpha_{k+1}\vert+\norm{\tilde\mu_{k+1}}),\beta_k+\alpha_k+1\big\}-\alpha_{k+1}\;\mbox{ as }\;k\ge 1.
\end{eqnarray*}
We clearly have that the sequence $\{\beta_k+\alpha_k\}$ is strictly increasing with $\beta_k+\alpha_k\ge k$ as $k\in\N$ and
\begin{eqnarray*}
\lim_{k\to\infty}\frac{\alpha_k}{\beta_k+\alpha_k}=\lim_{k\to\infty}\frac{\norm{\tilde\mu_k}}{\beta_k+\alpha_k}=0.
\end{eqnarray*}
Define $t_k:=(\beta_k+\alpha_k)^{-1}$ and $\mu_k:=(\tilde\mu_k+\beta_k\hat\mu)/(\beta_k+\alpha_k)\in K^\ast$ for all $k\in\N$ and observe that
$\frac{\mu_k-\hat\mu}{t_k}=\tilde\mu_k-\alpha_k\hat\mu$, which yields the relationships
\begin{eqnarray*}
\lim_{k\to\infty}\mu_k-\hat\mu=\lim_{k\to\infty}\frac{\tilde\mu_k-\alpha_k\hat\mu}{\beta_k+\alpha_k}=0\;\mbox{ and }\;\norm{w_k^\ast-\tilde w_k^\ast}\le k^{-1}\;\mbox{ with }\;\tilde w_k^\ast:=t_k^{-1}\nabla G(\zb)^\ast(\mu_k-\hat\mu).
\end{eqnarray*}
For each $t\in(0,1]$ we find $k\in\N$ such that $t\in(t_{k+1},t_k]$ and let $\mu_t:=((t_k-t)\mu_{k+1}+(t-t_{k+1})\mu_k)/(t_k-t_{k+1})$. This gives us $\mu_t\in K^\ast$, $\mu_t\to\hat\mu$ as $t\dn 0$, and
\begin{eqnarray*}
\nabla G(\zb)^\ast\frac{\mu_t-\hat\mu}{t}&=&\frac{(t_k-t)t_{k+1}}{(t_k-t_{k+1})t}\nabla G(\zb)^\ast\frac{\mu_{k+1}-\hat\mu}{t_{k+1}}+\frac{(t-t_{k+1})t_k}{(t_k-t_{k+1})t}\nabla G(\zb)^\ast\frac{\mu_k-\hat\mu}{t_k}\\
&=&\frac{(t_k-t)t_{k+1}}{(t_k-t_{k+1})t}\tilde w_{k+1}^\ast+\frac{(t-t_{k+1})t_k}{(t_k-t_{k+1})t}\tilde w_k^\ast.
\end{eqnarray*}
The nonnegativity of both $t_k-t$ and $t-t_{k+1}$  together with the relationships $(t_k-t)t_{k+1}+(t-t_{k+1})t_k=(t_k-t_{k+1})t$ and $w_k^\ast\to\tilde w_k^\ast$ as $k\to\infty$ imply the limiting conditions
\begin{equation}\label{EqAux_mu_t1}
\lim_{t\downarrow 0}\nabla_x G(\zb)^\ast\frac{\mu_t-\hat\mu}{t}=\lim_{t\downarrow 0}\pi_{X^\ast}\Big(\nabla G(\zb)^\ast\frac{\mu_t-\hat\mu}{t}\Big)=w_{x^\ast}^\ast,\quad\lim_{t\downarrow 0}\Big\la\nabla G(\zb)^\ast\frac{\mu_t-\hat\mu}{t},v\Big\ra=0.
\end{equation}
By $\skalp{z^\ast,v}=\skalp{\nabla G(\zb)^\ast\hat\mu,v}=0$, the latter condition yields
\begin{equation}\label{EqAux_mu_t2}
\skalp{\mu_t,\nabla G(\zb)v}=\skalp{\nabla G(\zb)^\ast\mu_t,v}=\skalp{\nabla G(\zb)^\ast\hat\mu,v}+\oo(t)=\oo(t).
\end{equation}

To proceed further, fix $\ve>0$ and consider the following optimization problem:
\begin{eqnarray*}
\min \norm{w}\;\mbox{ subject to }\;\la z^\ast,w\ra+\frac 12\Big\la \hat\mu,\nabla^2G(\zb)(v,v)\Big\ra\ge-\epsilon,\;{\rm dist}\Big(\nabla G(\zb)w+\frac 12\nabla^2G(\zb)(v,v);K\Big)\le\epsilon,
\end{eqnarray*}
which admits an optimal solution $w_\epsilon$ by Proposition~\ref{PropConicLinProg}. Note that the feasible region of this program is closed and convex. It is also easy to see that the function $\epsilon\to\norm{w_\epsilon}$ is decreasing, convex, and hence continuous. Considering first the case where $\norm{w_\epsilon}$ is not identically zero, for every $t>0$ define
\begin{eqnarray*}
\sigma_t:=\max\big\{t,\sup\{\norm{\nabla G(\zb)^\ast(\mu_\tau-\hat\mu)}\mv\tau\le t\}\big\}\;\mbox{ and }\;\epsilon_t:=\inf\big\{\epsilon>0\mv\sigma_t\norm{w_\epsilon}\le\epsilon\big\}.
\end{eqnarray*}
Since $\epsilon_t>0$ and $\norm{w_\epsilon}$ is continuous, we get
$\sigma_t\norm{w_{\epsilon_t}}=\epsilon_t$. Observe that the functions $t\to\sigma_t$ and $t\to\epsilon_t$ are increasing and that $\sigma_t\dn 0$ as $t\dn 0$, which allows us to claim that $\epsilon_t\dn 0$ as $t\dn 0$. Indeed, pick $\delta>0$ and $\bar t>0$ with $\sigma_{\bar t}<\delta/\norm{w_\delta}$. Then for every $0<t<\bar t$ we have $\sigma_t\norm{w_\delta}\le\sigma_{\bar t}\norm{w_\delta}<\delta$. This clearly yields $\epsilon_t<\delta$ and thus verifies the claim. Furthermore, it follows by setting $\bar w_t:=w_{\epsilon_t}$ that
\begin{subequations}
\begin{align}\label{EqProp_w_t_a}
&\liminf_{t\downarrow 0}\Big\la\hat\mu,\nabla G(\zb)\bar w_t+\frac 12\nabla^2G(\zb)(v,v)\Big\ra=\liminf_{t\downarrow 0}\la z^\ast,\bar w_t\ra+\frac 12\big\la\hat\mu,\nabla^2G(\zb)(v,v)\big\ra\ge 0,\\
\label{EqProp_w_t_c}&\lim_{t\downarrow 0}\big\la\mu_t-\hat\mu,\nabla G(\zb)\bar w_t\big\ra=\lim_{t\downarrow 0}\Big\la\mu_t-\hat\mu,\nabla G(\zb)\bar w_t+\frac 12\nabla^2G(\zb)(v,v)\Big\ra=0,\\
\nonumber&\lim_{t\downarrow 0}\Big\{{\rm dist}\Big(\nabla G(\zb)\bar w_t+\frac 12\nabla^2G(\zb)(v,v);K\Big)\Big\}=0
\end{align}
\end{subequations}
together with $t\norm{w_t}\to 0$ as $t\dn 0$. Observe also that these conditions hold with $\bar w_t=0$ if $\norm{w_\epsilon}=0$ for all $\epsilon>0$.
It follows from the latter condition above and from $\nabla G(\zb)v\in K$ that
\begin{eqnarray*}
{\rm dist}\big(G(\zb+tv+t^2\bar w_t;K\big)={\rm dist}\Big(t\nabla G(\zb)v+t^2\big(\nabla G(\zb)\bar w_t+\frac 12\nabla^2G(\zb)(v,v)\big)+\oo(t^2);K\Big)=\oo(t^2).
\end{eqnarray*}
Applying now MSCQ, for every $t>0$ we find $z_t=(p_t,x_t)$ with $G(z_t)\in K$ and $z_t=\zb+tv+t^2\bar w_t+\oo(t^2)$. Next choose $\alpha>0$ so small that $\alpha\norm{\nabla_{xx}^2\skalp{\hat\mu,G}(\zb)}<\frac 12$.
\if{, which ensures the estimate
\begin{eqnarray*}
\norm{q}^2+\alpha\nabla^2_{xx}\skalp{\hat\mu,G}(\zb)(q,q)\ge\frac 12\norm{q}^2\;\mbox{ for all }\;q\in X.
\end{eqnarray*}}\fi
Since $\pi_{X^\ast}(z^\ast)=\xba$ and $z^\ast=\nabla G(\zb)^\ast\hat\mu$, we have $\nabla_xG(\zb)^\ast\hat\mu=\xba$ and deduce from the second-order sufficient condition in \cite[Theorem~3.63]{BonSh00} that $\xb$ is a strict local solution to the optimization problem
\begin{equation}\label{EqAux_QP}
\min\frac 12\norm{\xb+\alpha\J_X(\xba)-x}^2\;\mbox{ subject to }\;G(\pb,x)\in K,
\end{equation}
where $\J_X\colon X^\ast\to X$ stands for the classical Riesz isomorphism. Choose $\rho>0$ so that $\xb$ is the unique global solution to this program on $\B_\rho(\xb)$ and for every $t\in(0,1]$ denote by $\xb_t$ a global solution to
\begin{equation}\label{EqAuxQP_t}
\min\frac 12\norm{x_t+\alpha\J_X(x^\ast_t)-x}^2\;\mbox{ subject to }\;x\in\Gamma(p_t)\cap\B_\rho(\xb)
\end{equation}
with $x^\ast_t:=\nabla_xG(z_t)^\ast\mu_t$. Observe that $x_t$ is feasible to the latter program for all $t>0$ sufficiently small. It is not hard to check to that $\xb_t$ converges to $\xb$ as $t\downarrow 0$. Indeed, suppose that $\norm{\xb_{t_k}-\xb}\ge\delta>0$ for $t_k\downarrow 0$, and so $\xb_{t_k}$ converges to some $\tilde x$ along a subsequence with $G(\pb,\tilde x)=\lim_{k\to\infty}G(p_{t_k},\xb_{t_k})\in K$ and $\delta\le\norm{\tilde x-\xb}\le\rho$. Then passing to the limit in the relationships
\begin{eqnarray*}
\frac 12\norm{x_{t_k}+\alpha\J_X(x^\ast_{t_k})-\xb_{t_k}}^2\le\frac 12\norm{x_{t_k}+\alpha\J_X(x^\ast_{t_k})-x_{t_k}}^2=\frac{\alpha^2}2\norm{\J_X(x^\ast_{t_k})}^2
\end{eqnarray*}
gives us the inequality $\frac 12\norm{\xb+\J_X(\xba)-\tilde x}^2\le\frac 12\norm{\xb+\J_X(\xba)-\xb}^2$, which contradicts the uniqueness of the minimizer $\xb$ for \eqref{EqAux_QP} on $B_\rho(\xb)$. Hence we get $\xb_t\to\xb$ and $\norm{\xb_t-x_t}\to 0$ as $t\dn 0$.

The feasibility of $x_t$ in \eqref{EqAuxQP_t} yields $\norm{x_t+\alpha \J_X(x^\ast_t)-\xb_t}^2\le\norm{\alpha\J_X(x^\ast_t)}^2$, and thus
\begin{eqnarray}\label{EqAux_Diff_x_t}
&&\norm{x_t-\xb_t}^2\le 2\alpha\skalp{x^\ast_t,\xb_t-x_t}=2\alpha\skalp{\mu_t,\nabla_xG(z_t)(\xb_t-x_t)}\\
\nonumber&&=2\alpha\skalp{\mu_t,G(p_t,\xb_t)-G(p_t,x_t)}-\alpha\nabla_{xx}^2\skalp{\mu_t,G}(p_t,x_t)(\xb_t-x_t,\xb_t-x_t)+\oo(\norm{\xb_t-x_t}^2).
\end{eqnarray}
It follows from \eqref{EqAux_mu_t2}, \eqref{EqProp_w_t_a}, \eqref{EqProp_w_t_c}, and $\lim_{t\downarrow 0}\mu_t=\hat\mu$ that
\begin{eqnarray*}
\big\la\mu_t,G(p_t,x_t)\big\ra&=&t\big\la\mu_t,\nabla G(\zb)v\big\ra+t^2\Big\la\mu_t,\nabla G(\zb)\bar w_t+\frac 12\nabla^2G(\zb)(v,v)\Big\ra+\oo(t^2)\\
&=&t\big\la\mu_t,\nabla G(\zb)v\big\ra+t^2\Big\la\hat\mu,\nabla G(\zb)\bar w_t+\frac 12\nabla^2G(\zb)(v,v)\Big\ra+\oo(t^2)\geq \oo(t^2).
\end{eqnarray*}
On the other hand, we have $\la\mu_t,G(p_t,\xb_t)\ra\le 0$ due to $G(p_t,\xb_t)\in K$ and $\mu_t\in K^\ast$. Combining this with \eqref{EqAux_Diff_x_t} verifies the validity of the estimate
\begin{eqnarray*}
\norm{x_t-\xb_t}^2\le-\alpha\nabla_{xx}^2\skalp{\mu_t,G}(p_t,x_t)(\xb_t-x_t,\xb_t-x_t)+\oo(\norm{\xb_t-x_t}^2)+\oo(t^2).
\end{eqnarray*}
Since $\lim_{t\downarrow 0}\nabla_{xx}^2\skalp{\mu_t,G}(p_t,x_t)=\nabla_{xx}^2\skalp{\hat\mu,G}(\zb)$, we get $\alpha \norm{\nabla_{xx}^2\skalp{\mu_t,G}(p_t,x_t)}\le\frac 58$ whenever $t>0$ is small, and hence $-\alpha\nabla_{xx}^2\skalp{\mu_t,G}(p_t,x_t)(\xb_t-x_t,\xb_t-x_t)\le\frac 58\norm{x_t-\xb_t}^2$. Then $\frac 14\norm{x_t-\xb_t}^2\le\oo(t^2)$ and
\begin{equation}\label{EqAux_Diff_x_t2}
\norm{x_t-\xb_t}=\oo(t)\;\mbox{ and }\;\norm{(\xb_t,p_t)-(\zb+tv)}=\oo(t)
\end{equation}
for all $t>0$ sufficiently small. Furthermore, it follows from \eqref{EqAux_mu_t1} that
\begin{eqnarray*}
x^\ast_t&=&\nabla_xG(z_t)^\ast\mu_t=\nabla_xG(\zb)^\ast\hat\mu+\big(\nabla_xG(z_t)^\ast\hat\mu-\nabla_xG(\zb)^*\hat\mu\big)+\nabla_xG(z_t)^\ast(\mu_t-\hat\mu)\\
&=&\xba+t\Big(\nabla\big(\nabla_xG(\cdot)^\ast\hat\mu\big)(\zb)v+\nabla_xG(\zb)^\ast\frac{\mu_t-\hat\mu}t\Big)+\oo(t)\\
&=&\xba+t\big(\nabla(\nabla_xG(\cdot)^\ast\hat\mu)(\zb)v+w_{x^\ast}^\ast\big)+\oo(t).
\end{eqnarray*}
Remembering that $\ox_t\to\ox$ as $t\dn 0$ and using the well-known necessary optimality condition (see, e.g., \cite[Theorem 6.12]{RoWe98}) for the optimal solution $\ox_t$ to program \eqref{EqAuxQP_t} tell us that
\begin{eqnarray*}
0\in\J_X^{-1}(\xb_t-x_t)-\alpha x^\ast_t+N_{\Gamma(p_t)}(\xb_t)\;\mbox{ for all small }\;t>0.
\end{eqnarray*}
Since $N_{\Gamma(p_t)}(\xb_t)$ is a cone, we obtain from the above and the left relationship in \eqref{EqAux_Diff_x_t2} that
\begin{eqnarray*}
{\rm dist}\big(\xba+t\big(\nabla(\nabla_xG(\cdot)^\ast\hat\mu)(\zb)v+w_{x^\ast}^\ast\big);N_{\Gamma(p_t)}(\xb_t)\big)={\rm dist}\big(x^\ast_t; N_{\Gamma(p_t)}(\xb_t)\big)+\oo(t)=\oo(t),
\end{eqnarray*}
which being combined with the right relationship therein verifies the first inclusion of the theorem with the tangent derivability. To get finally the refined inclusion \eqref{RemConv_w_k}, it remains to apply the first assertion of the theorem by setting $w_k^\ast=w^\ast$ for all $k$ together with the equality $N_{\K_{{\rm gph}\,\Gamma}(\zb,z^\ast)}(v)=\K_{{\rm gph}\,\Gamma}(\zb,z^\ast)^\ast\cap[v]^\perp$, which is due to the duality in \eqref{tan1}. $\h$\vspace*{-0.15in}

\section{Graphical Derivative of the Normal Cone via Reduction Data}\sce\vspace*{-0.1in}

In this section we derive upper estimates and exact expressions for the graphical derivative of the normal cone mapping \eqref{EqNormalConeMap} in terms of the auxiliary data of the parametric reduced system $G(z)\in K$ with $z=(p,x)$. The Robinson stability \eqref{EqRS} of the latter systems plays a crucial role in our consideration.

For the subsequent analysis we need the following {\em characteristic subspace} of the parameter space $P$ defined at the reference point $\oz=(\op,\ox)$ by
\begin{eqnarray}\label{char}
\Pp:=\big\{q\in P\big|\;\nabla_p G(\zb)q\in\rge\nabla_xG(\zb)+\lin K\big\},
\end{eqnarray}
where $\rge A$ stands for the range of the linear operator. The next lemma on descriptions of tangent vectors to the graph of \eqref{EqNormalConeMap}
is basic for further developments. To avoid any confusion, note that the notation $N_{T_{{\rm gph}\,\Gamma}(\zb)}(v)$ and similar ones in what follows mean that we consider the normal cone at the point $v$ to the tangent cone to the graph of $\Gamma$ taken at the point $\oz$.\vspace*{-0.07in}

\begin{lemma}[\bf descriptions of graphical tangent directions]\label{PropTanDir} Assume that the reduced system $G(p,x)\in K$ enjoys the Robinson stability property at $(\pb,\xb)$. Then for every graphical tangent direction $(q,u,u^\ast)\in T_{{\rm gph}\,\Psi}(\pb,\xb,\xba)$ the following assertions hold:

{\bf(i)} There exists $\hat\mu\in N_K(\nabla G(\zb)v)$ with $v:=(q,u)$ such that $\nabla_xG(\zb)^\ast\hat\mu=\xba$ and
\begin{equation}\label{EqBasicPropGraphDer}
u^\ast\in\nabla\big(\nabla_xG(\cdot)^\ast\hat\mu\big)(\zb)v+\K_{\Gamma(\pb)}(\xb,\xba)^*.
\end{equation}

{\bf(ii)} If either $\Pp=P$ or the cone $K$ is polyhedral, then there exist $z^\ast\in N_{T_{{\rm gph}\,\Gamma}(\zb)}(v)\cap\pi_{X^\ast}^{-1}(\xba)$ with $\pi_{X^*}$ taken from \eqref{proj} and a directional multiplier $\hat\mu\in\M(\zb,z^\ast;v)$ from \eqref{mult} such that
\begin{equation}\label{EqDesiredPropGraphDer}
u^\ast\in\nabla\big(\nabla_xG(\cdot)^\ast\hat\mu\big)(\zb)v+\pi_{X^\ast}\big(N_{\K_{{\rm gph}\,\Gamma}(\zb,z^\ast)}(v)\big).
\end{equation}

{\bf(iii)} If there exist sequences $t_k\downarrow 0$, $(q_k,u_k,u_k^\ast)\to(q,u,u^\ast)$ with $\xba+t_ku_k^\ast\in N_{\Gamma(\pb+t_kq_k)}(\xb+t_ku_k)$ and ${\rm dist}(q_k;q+\Pp)=\OO(t_k)$, then we can find $z^\ast\in N_{T_{{\rm gph}\,\Gamma}(\zb)}(v)\cap\pi_{X^\ast}^{-1}(\xba)$, $\hat\mu\in\M(\zb,z^\ast;v)$, and
$w_{x^\ast}^\ast\in X^\ast$ for which there is a sequence of $w_k^\ast\in\K_{{\rm gph}\,\Gamma}(\zb,z^\ast)^\ast$ such that
\begin{eqnarray*}
w_{x^\ast}^\ast=\lim_{k\to\infty}\pi_{X^\ast}(w_k^\ast),\;\lim_{k\to\infty}\skalp{w_k^\ast,v}=0,\;\mbox{ and }\;u^\ast=\nabla\big(\nabla_xG(\cdot)^\ast\hat\mu\big)(\zb)v+w_{x^\ast}^\ast.
\end{eqnarray*}
\end{lemma}\vspace*{-0.07in}
{\bf Proof.} To verify assertion (i), take $t_k\downarrow 0$, $(q_k,u_k,u_k^\ast)\to(q,u,u^\ast)$ with $\xba+t_ku_k^\ast\in N_{\Gamma(\pb+t_kq_k)}(\xb+t_ku_k)$ for all $k\in\N$. The robustness of the Robinson stability property ensures that the system $0\in G(\pb+t_kq_k,\cdot)-K$ fulfills MSCQ with the uniform modulus $\kappa$ at $\xb+t_ku_k$ for all $k$ sufficiently large. It follows from Proposition~\ref{LemNormalCone} that there are
$\mu_k\in N_K(G(\zb+t_kv_k))\cap\kappa\norm{\xba+t_ku_k^\ast}\B_{E^\ast}$ satisfying $\xba+t_ku_k^\ast=\nabla_xG(\zb+t_kv_k)^\ast\mu_k$ with $v_k:=(q_k,u_k)$. Since the set $\{\xba+t_ku_k^\ast|\;k\in\N\}$ is clearly bounded, so is $\{\mu_k|\;k\in\N\}$, which therefore contains a convergent subsequence $\mu_k\to\hat\mu$ to some multiplier $\hat\mu$ satisfying $\hat\mu\in N_K(G(\zb))=K^\ast$ and $\nabla_xG(\zb)^\ast\hat\mu=\xba$. By taking into account that $G(\xb+t_kv_k)=t_k\nabla G(\zb)v_k+\oo(t_k)\in K$ and the conic structure of $K$, we get $\nabla G(\zb)v_k+\oo(t_k)/t_k\in K$ and hence $\nabla G(\zb)v\in K$. Furthermore, $N_K(G(\zb+t_kv_k))=K^\ast\cap[G(\zb+t_kv_k)]^\perp$, which yields $\skalp{\mu_k,G(\zb+t_kv_k)}=0$ for all $k$ and therefore
\begin{eqnarray*}
0=\lim_{k\to\infty}\frac{\big\la\mu_k,G(\zb+t_kv_k)\big\ra}{t_k}=\big\la\hat\mu,\nabla G(\zb)v\big\ra.
\end{eqnarray*}
This gives us $\hat\mu\in[\nabla G(\zb)v]^\perp$ verifying $\hat\mu\in K^\ast\cap[\nabla G(\zb)v]^\perp=N_K(\nabla G(\zb)v)$. Observe also that
\begin{eqnarray*}
u_k^\ast&=&\frac{\nabla_xG(\zb+t_kv_k)^\ast\mu_k-\xba}{t_k}=\frac{\nabla_xG(\zb+t_kv_k)^\ast\mu_k-\nabla_xG(\zb)^\ast\hat\mu}{t_k}\\
&=&\nabla\big(\nabla_xG(\cdot)^\ast\hat\mu\big)(\zb)v+\frac{\oo(t_k)}{t_k}+\nabla_xG(\zb)^\ast\frac{\mu_k-\hat\mu}{t_k},
\end{eqnarray*}
which justifies the validity of the limiting representation
\begin{equation}\label{EqAuxBasicProp}
u^\ast-\nabla\big(\nabla_xG(\cdot)^\ast\hat\mu\big)(\zb)v=\lim_{k\to\infty}\nabla_xG(\zb)^\ast\frac{\mu_k-\hat\mu}{t_k}.
\end{equation}
Picking now an arbitrary critical direction $w\in\K_{\Gamma(\pb)}(\xb,\xba)=T_{\Gamma(\pb)}(\xb)\cap[\xba]^\perp$, we get $\nabla G_x(\zb)w\in K$ and $\skalp{\xba,w}=\skalp{\nabla_xG(\zb)^\ast\hat\mu,w}=0$. Together with $\mu_k\in K^\ast$ it implies that
\begin{eqnarray*}
0\ge\frac{\skalp{\mu_k,\nabla _xG(\zb)w}}{t_k}=\frac{\skalp{\nabla_xG(\zb)^\ast\mu_k,w}}{t_k}=\Big\la\nabla_xG(\zb)^\ast\frac{\mu_k-\hat\mu}{t_k},w\Big\ra
\end{eqnarray*}
and consequently verifies that $\nabla_xG(\zb)^\ast\frac{\mu_k-\hat\mu}{t_k}\in\K_{\Gamma(\pb)}(\xb,\xba)^\ast$. This deduce \eqref{EqBasicPropGraphDer} from \eqref{EqAuxBasicProp} due to the polar cone closedness thus completing the proof of assertion (i).

For the rest of the proof of this theorem we set $z^\ast:=\nabla G(\zb)^\ast\hat\mu\in N_{{\rm gph}\,\Gamma}(\zb)\cap\pi_{X^\ast}^{-1}(\xba)$. As already shown, $\skalp{z^\ast,v}=\skalp{\hat\mu,\nabla G(\zb)v}=0$ and hence $z^\ast\in N_{{\rm gph}\,\Gamma}(\zb)\cap[v]^\perp\cap\pi_{X^\ast}^{-1}(\xba)=N_{T_{{\rm gph}\, \Gamma}(\zb)}(v)\cap\pi_{X^\ast}^{-1}(\xba)$. We also have $v\in T_{{\rm gph}\,\Gamma}(\zb)\cap[z^\ast]^\perp=\K_{{\rm gph}\,\Gamma}(\zb,z^\ast)$ by definition \eqref{cc}.

To proceed next with justifying assertion (ii), let us first verify that the sequence $\big\{q_k^\ast:=\nabla_pG(\zb)^\ast\frac{\mu_k-\hat \mu}{t_k}\big\}$ is bounded. Starting with the case where $\Pp=P$ and arguing by contradiction, suppose that $\{q_k^\ast\}$ is unbounded. Then passing to a subsequence if necessary gives us $\tilde q\in P$ such that
\begin{eqnarray*}
\big\la q_k^\ast,\tilde q\big\ra=\Big\la\frac{\mu_k-\hat\mu}{t_k},\nabla_pG(\zb)\tilde q\Big\ra\to\infty\;\mbox{ as }\;k\to\infty.
\end{eqnarray*}
The assumption of $\Pp=P$ and the construction of $\Pp$ in \eqref{char} allow us to find $\tilde u\in X$ and $\tilde\mu\in\lin K$ with $\nabla_pG(\zb)\tilde q=\nabla_xG(\zb)\tilde u+\tilde\mu$. It follows from $\mu_k-\hat\mu\in K^\ast-K^\ast=(\lin K)^\perp$ and \eqref{EqAuxBasicProp} that
\begin{eqnarray*}
\infty&=&\lim_{k\to\infty}\Big\la\frac{\mu_k-\hat\mu}{t_k},\nabla_pG(\zb)\tilde q\Big\ra=\lim_{k\to\infty}\Big\la\frac{\mu_k-\hat\mu}{t_k},\nabla_xG(\zb)\tilde u+\tilde\mu\Big\ra=\lim_{k\to\infty}\Big\la\nabla_xG(\zb)^\ast\frac{\mu_k-\hat\mu}{t_k},\tilde u\Big\ra\\
&=&\big\la u^\ast-\nabla\big(\nabla_xG(\cdot)^\ast\hat\mu\big)(\zb)v,\tilde u\big\ra,
\end{eqnarray*}
a contradiction that verifies the boundedness of $\{q_k^\ast\}$ if $\Pp=P$. It remains to consider the case where $K$ is polyhedral. Then its polar cone $K^\ast$ is polyhedral as well and the normal cones $N_K(G(\zb+t_kv_k))=K^\ast\cap[G(\zb+t_kv_k)]^\perp$  for each $k$ are faces of $K^\ast$. Since a polyhedral cone has only finitely many faces, by passing to a subsequence we may suppose that all these faces are identical, i.e., there is $s\in K$ such that $N_K(G(\zb+t_kv_k))=K^\ast\cap[s]^\perp$ for all $k\in\N$. The classical result by Walkup and Wets \cite{WalkWe69} tells us that the set-valued mapping $b^\ast\rightrightarrows{\cal S}(b^\ast):=\{\mu\in K^\ast\cap[s]^\perp|\,\nabla_xG(\zb)^\ast\mu=b^\ast\}$ is Lipschitz continuous in the Hausdorff metric provided that it is nonempty. Since $\mu_k\in{\cal S}(\nabla_xG(\zb)^\ast\mu_k)$ and $\hat\mu\in{\cal S}(\xba)$, there is a constant $L>0$ such that for every $k\in\N$ we can find $\tilde\mu_k\in{\cal S}(\nabla_xG(\zb)^\ast\mu_k)$ for which
\begin{eqnarray*}
\norm{\tilde\mu_k-\hat\mu}&\le&L\norm{\nabla_xG(\zb)^\ast\mu_k-\xba}=L\norm{\xba+t_ku_k^\ast+\big(\nabla_xG(\zb)^\ast-\nabla_xG(\zb+t_kv_k)^\ast\big)\mu_k-\xba}\\
&=&L\big(t_k\norm{u_k^\ast-\nabla\big(\nabla_xG(\cdot)^\ast\mu_k\big)(\zb)v_k}+\oo(t_k)\big),
\end{eqnarray*}
and so $\{(\tilde\mu_k-\hat\mu)/t_k\}$ is bounded. Defining further $\tilde u_k^\ast:=(\nabla_xG(\zb+t_kv_k)^\ast\tilde\mu_k-\xba)/t_k$, we get
\begin{eqnarray*}
\xba+t_k\tilde u_k^\ast=\nabla_xG(\zb+t_kv_k)^\ast\tilde\mu_k\in N_{\Gamma(\pb+t_kq_k)}(\xb+t_ku_k)\;\mbox{ and}
\end{eqnarray*}\vspace*{-0.4in}
\begin{eqnarray*}
\lim_{k\to\infty}(\tilde u_k^\ast-u_k^\ast)&=&\lim_{k\to\infty}\nabla_xG(\zb+t_kv_k)^\ast\frac{\tilde\mu_k-\mu_k}{t_k}\\
&=&\lim_{k\to\infty}\nabla_xG(\zb)^\ast\frac{\tilde\mu_k-\mu_k}{t_k}+\Big(\nabla\big(\nabla_xG(\zb)^\ast\big)v_k+\frac{\oo(t_k)}{t_k}\Big)(\tilde\mu_k-\mu_k)=0.
\end{eqnarray*}
Hence $\tilde u_k^\ast\to u^\ast$ and we can rename the sequences $\{\tilde u_k^\ast\}$ and $\{\tilde\mu_k\}$ by $\{u_k^\ast\}$ and $\{\mu_k\}$, respectively. Then the boundedness of $\big\{\frac{\tilde\mu_k-\hat\mu}{t_k}\big\}$ yields the one for $\{q_k^\ast\}$ in the polyhedral case.

The boundedness of the sequences $\big\{q_k^\ast=\nabla_pG(\zb)^\ast\frac{\mu_k-\hat\mu}{t_k}\big\}$ and $\big\{\nabla_xG(\zb)^\ast\frac{\mu_k-\hat\mu}{t_k}\big\}$ implies this property for the sequence $\big\{\nabla G(\zb)^\ast\frac{\mu_k-\hat\mu}{t_k}\big\}$, and thus we have $\nabla G(\zb)^\ast\frac{\mu_k-\hat\mu}{t_k}\to w^\ast$ as $k\to\infty$ along a subsequence. Fix $w\in\K_{{\rm gph}\,\Gamma}(\zb,z^\ast)$ and get that
\begin{eqnarray*}
0\ge\frac{\big\la\mu_k,\nabla G(\zb)w\big\ra}{t_k}=\Big\la\nabla G(\zb)^\ast\frac{\mu_k-\hat\mu}{t_k},w\Big\ra\;\mbox{ for any }\;k\in\N,
\end{eqnarray*}
which yields $\skalp{w^\ast,w}\le 0$ and $w^\ast\in\K_{{\rm gph}\,\Gamma}(\zb,z^\ast)^\ast$. On the other hand, we have
\begin{eqnarray*}
0&\le&\liminf_{k\to\infty}\frac{\big\la\mu_k-\hat\mu,G(\zb+t_kv_k\big\ra}{t_k^2}\\
&=&\liminf_{k\to\infty}\left(\Big\la\frac{\mu_k-\hat\mu}{t_k},\nabla G(\zb)v_k\Big\ra+\frac 12 \Big\la\mu_k-\hat\mu,\nabla^2G(\zb)(v_k,v_k)+\frac{\oo(t_k^2)}{t_k^2}\Big\ra\right)=\skalp{w^\ast,v}
\end{eqnarray*}
verifying that $w^\ast\in\K_{{\rm gph}\,\Gamma}(\zb,z^\ast)^\ast\cap[v]^\perp=N_{\K_{{\rm gph}\,\Gamma}(\zb,z^\ast)}(v)$ by the duality correspondence \eqref{tan1}. This shows that condition \eqref{EqDesiredPropGraphDer} follows from \eqref{EqAuxBasicProp}. To justify (ii), it remains to check the inclusion $\hat\mu\in\M(\zb,z^\ast;v)$. Indeed, for every $\mu\in\M(\zb,z^\ast)$ we have $\nabla G(\zb)^\ast\mu=\nabla G(\zb)^\ast\hat\mu$ and therefore
\begin{eqnarray*}
0&\le&\liminf_{k\to\infty}t_k^{-2}\big\la\mu_k-\mu,G(\zb+t_kv_k)\big\ra\\
&=&\liminf_{k\to\infty}\left(\Big\la\nabla G(\zb)^\ast\frac{\mu_k-\mu}{t_k},v_k\Big\ra+\frac 12\Big\la\mu_k-\mu,\nabla^2G(\zb)(v_k,v_k)+\frac{\oo(t_k^2)}{t_k^2}\Big\ra\right)\\
&=&\liminf_{k\to\infty}\left(\Big\la\nabla G(\zb)^\ast\frac{\mu_k-\hat\mu}{t_k},v_k\Big\ra+\frac 12\Big\la\mu_k-\mu,\nabla^2G(\zb)(v_k,v_k)+\frac{\oo(t_k^2)}{t_k^2}\Big\ra\right)\\
&=&\skalp{w^\ast,v}+\frac 12\big\la\hat\mu-\mu,\nabla^2G(\zb)(v,v)\big\ra=\frac 12\big\la\hat\mu-\mu,\nabla^2G(\zb)(v,v)\big\ra,
\end{eqnarray*}
which justifies the inclusion $\hat\mu\in\M(\zb,z^\ast;v)$ and completes the proof of assertion (ii).

Now we turn to the final assertion (iii) of the theorem. Represent $P$ as the direct sum $P=P\bigoplus\hat\Pp$ with some subspace $\hat\Pp$ and denote by $\pi_{\hat \Pp}$ the corresponding projection of $P$ onto $\hat\Pp$ along $\Pp$. For each $k\in\N$ take $\tilde q_k\in \Pp$ with $\norm{q+\tilde q_k-q_k}={\rm dist}(q_k;q+\Pp)$ and get $\pi_{\hat\Pp}(q-q_k)=\pi_{\hat\Pp}(q+\tilde q_k-q_k)$ by the construction. The continuity of the projection and the distance assumption in (iii) yields the representation $\hat q_k:=\pi_{\hat\Pp}(q-q_k)=\OO(t_k)$. We claim that
\begin{eqnarray}\label{cla}
\lim_{k\to\infty}\Big\la\nabla_pG(\zb)^\ast\frac{\mu_k-\hat\mu}{t_k},q-q_k-\hat q_k\Big\ra=0.
\end{eqnarray}
Indeed, the failure of \eqref{cla} together with $q-q_k-\hat q_k\in\Pp$ and $q-q_k-\hat q_k\to 0$ gives us $\tilde q\in\Pp$ with
\begin{eqnarray*}
\Big\vert\Big\la\nabla_pG(\zb)^\ast\frac{\mu_k-\hat\mu}{t_k},\tilde q\Big\ra\Big\vert=\Big\vert\Big\la\frac{\mu_k-\hat\mu}{t_k},\nabla_pG(\zb)\tilde q\Big\ra\Big\vert\to\infty\;\mbox{ as }\;k\to\infty.
\end{eqnarray*}
By $\tilde q\in\Pp$ we find $\tilde u\in X$ and $\tilde \mu\in \lin K$ satisfying $\nabla_pG(\zb)\tilde q=\nabla_xG(\zb)\tilde u+\tilde\mu$ and then proceed as the proof of (ii) to get the contradiction that  verifies \eqref{cla}.

Next define $w_k^\ast:=\nabla G(\zb)^\ast\frac{\mu_k-\hat\mu}{t_k}$ for $k\in\N$ and show that $\{w_k^\ast\}$ is the sequence whose existence is claimed in (iii). Observe that $\skalp{w_k^\ast,w}=\skalp{\nabla G(\zb)^\ast\frac{\mu_k-\hat\mu}{t_k},w}=\frac{\skalp{\nabla G(\zb)^\ast\mu_k,w}}{t_k}\le 0$ for every $w\in\K_{{\rm gph}\,\Gamma}(\zb,z^\ast)$, and so $w_k^\ast\in\K_{{\rm gph}\,\Gamma}(\zb,z^\ast)^\ast$. It follows from the boundedness of $\big\{\nabla_xG(\zb)^\ast\frac{\mu_k-\hat\mu}{t_k}\big\}$ and $\big\{\frac{\hat q_k}{t_k}\big\}$ that
\begin{eqnarray*}
0&\le&\liminf_{k\to\infty}t_k^{-2}\big\la\mu_k-\hat\mu,G(\zb+t_kv_k)\big\ra\\
&=&\liminf_{k\to\infty}\left(\Big\la\frac{\mu_k-\hat\mu}{t_k},\nabla G(\zb)v_k\Big\ra+\frac12\Big\la\mu_k-\hat\mu,\nabla^2G(\zb)(v_k,v_k)+\frac{\oo(t_k)^2}{t_k^2}\Big\ra\right)\\
&=&\liminf_{k\to\infty}\Big\la\frac{\mu_k-\hat\mu}{t_k},\nabla G(\zb)v_k\Big\ra\\
&=&\liminf_{k\to\infty}\left(\Big\la\nabla G(\zb)^\ast\frac{\mu_k-\hat\mu}{t_k},v\Big\ra-\Big\la\nabla_xG(\zb)^\ast\frac{\mu_k-\hat\mu}{t_k},u-u_k\Big\ra\right.\\
&&\qquad\qquad\left.-\Big\la\nabla_pG(\zb)^\ast\frac{\mu_k-\hat\mu}{t_k},q-q_k-\hat q_k\Big\ra-\Big\la\nabla_pG(\zb)^\ast(\mu_k-\hat\mu),\frac{\hat q_k}{t_k}\Big\ra\right)\\
&=&\liminf_{k\to\infty}\skalp{w_k^\ast,v}\le\limsup_{k\to\infty}\skalp{w_k^\ast,v}\le 0,
\end{eqnarray*}
which implies the limiting condition $\lim_{k\to\infty}\skalp{w_k^\ast,v}=0$ asserted in (iii) and $\liminf_{k\to\infty}\big\la\frac{\mu_k-\hat\mu}{t_k},\nabla G(\zb)v_k\big\ra=0$. To show that $\hat\mu\in\M(\zb,z^\ast;v)$, take any $\mu\in\M(\zb,z^\ast)$ and use $\nabla G(\zb)^\ast\mu=\nabla G(\zb)^\ast\hat\mu$ to get
\begin{eqnarray*}
0&\le&\liminf_{k\to\infty}t_k^{-2}\big\la\mu_k-\mu,G(\zb+t_kv_k)\big\ra\\
&=&\liminf_{k\to\infty}\left(\Big\la\frac{\mu_k-\mu}{t_k},\nabla G(\zb)v_k\Big\ra+\frac 12\Big\la\mu_k-\mu,\nabla^2G(\zb)(v_k,v_k)+\frac{\oo(t_k)^2}{t_k^2}\Big\ra\right)\\
&=&\liminf_{k\to\infty}\left(\Big\la\frac{\mu_k-\hat\mu}{t_k},\nabla G(\zb)v_k\Big\ra+\frac 12\Big\la\mu_k-\mu,\nabla^2G(\zb)(v_k,v_k)+\frac{\oo(t_k)^2}{t_k^2}\Big\ra\right)\\
&=&\frac 12\big\la\hat\mu-\mu,\nabla^2G(\zb)(v,v)\big\ra,
\end{eqnarray*}
which yields the desired inclusion $\hat\mu\in\M(\zb,z^\ast;v)$ and thus completes the proof of the lemma.$\h$

Now we present the main result of this section on calculating the graphical derivative of the normal cone mapping \eqref{EqNormalConeMap} via the reduction data that is a direct consequence of those obtained above.\vspace*{-0.03in}

\begin{Theorem}[\bf computing the graphical derivative of the normal cone mapping for PVS via the reduction data]\label{ThMainTh} In addition to the standing assumption, suppose that the reduced system $G(p,x)\in K$ enjoys RS at $(\pb,\xb)$ and that either $\Pp=P$ or $K$ is polyhedral. Then we have
\begin{eqnarray*}
T_{{\rm gph}\,\Psi}(\pb,\xb,\xba)&=&\Big\{(q,u,u^\ast)\Big|\;\exists\,z^\ast\in N_{T_{{\rm gph}\,\Gamma}(\zb)}(q,u)\cap\pi_{X^\ast}^{-1}(\xba),\;\hat\mu\in\M\big(\zb,z^\ast;(q,u)\big)\\
&&\quad\mbox{with }\;u^\ast\in\nabla\big(\nabla_x G(\cdot)^\ast\hat\mu\big)(\zb)(q,u)+\pi_{X^\ast}\big(N_{\K_{{\rm gph}\,\Gamma}(\zb,z^\ast)}(q,u)\big)\Big\}.
\end{eqnarray*}
Consequently, for all $v=(q,u)\in Z$ we have the graphical derivative formula
\begin{eqnarray*}
D\Psi(\pb,\xb,\xba)(v)&=&\Big\{\nabla\big(\nabla_x G(\cdot)^\ast\hat\mu\big)(\zb)v+\pi_{X^\ast}\big(N_{\K_{{\rm gph}\,\Gamma}(\zb,z^\ast)}(v)\big)\Big|\\
&&\quad z^\ast\in N_{T_{{\rm gph}\,\Gamma}(\zb)}(v)\cap\pi_{X^\ast}^{-1}(\xba),\;\hat\mu\in\M(\zb,z^\ast;v)\Big\}.
\end{eqnarray*}
\end{Theorem}\vspace*{-0.05in}
{\bf Proof.} It follows from inclusion \eqref{RemConv_w_k} in Theorem~\ref{PropSuffTanDir} and the results of Lemma~\ref{PropTanDir}(ii), by taking into account that RS for the system $G(p,x)\in K$ at $(\pb,\xb)$ implies its MSCQ in Theorem~\ref{PropSuffTanDir}. $\h$\vspace*{0.03in}

To conclude this section, we derive an upper estimate and a precise formula for computing $D\Psi$ via the reduced data that does not impose any polyhedrality assumption on $K$ while replacing it by {\em strict complementarity}. The next lemma comes first.\vspace*{-0.03in}

\begin{lemma}[\bf critical cone under strict complementarity]\label{LemStrictCompl}
Assume that both MSCQ and strict complementarity conditions are satisfied at $\xb$ and $(\xb,\xba)$, respectively. Then for all $v=(q,u)\in\Pp\times X$ and $\hat\mu\in N_K(\nabla G(\zb)v)$ with $\nabla_xG(\zb)^\ast\hat\mu=\xba$ we have the inclusion
\begin{eqnarray*}
\K_{\Gamma(\pb)}(\xb,\xba)^\ast\subset\pi_{X^\ast}\big(N_{\K_{{\rm gph}\,\Gamma}(\zb,\nabla G(\zb)^\ast\hat\mu)}(v)\big).
\end{eqnarray*}
\end{lemma}\vspace*{-0.03in}
{\bf Proof.} Pick $w_{x^\ast}\in\K_{\Gamma(\pb)}(\xb,\xba)^\ast$ and get by Proposition~\ref{LemNormalCone} that
$\K_{\Gamma(\pb)}(\xb,\xba)^\ast=\nabla_xG(\zb)^\ast K^\ast+[\xba]$. This gives us $\mu\in K^\ast$ and $\alpha\in\R$ such that $w_{x^\ast}^\ast=\nabla_xG(\zb)^\ast\mu+\alpha\xba=\nabla_xG(\zb)^\ast(\mu+\alpha\hat\mu)$. By $q\in\Pp$ we get $\hat u\in X$ and $\th\in\lin K$ satisfying $\nabla_pG(\zb)q=\nabla_xG(\zb)\hat u+\th$, which yields
\begin{eqnarray*}
\nabla G(\zb)v=\nabla_xG(\zb)u+\nabla_pG(\zb)q=\nabla_xG(\zb)(u+\hat u)+\th\in K\;\mbox{ and }\;\nabla_xG(\zb)(u+\hat u)\in K.
\end{eqnarray*}
It follows from $\skalp{\hat\mu,\th}=0$ and $\hat\mu\in N_K(\nabla G(\zb)v)$ that
\begin{eqnarray*}
0=\big\la\hat\mu,\nabla G(\zb)v\big\ra=\big\la\hat\mu,\nabla_xG(\zb)(u+\hat u)+\th\big\ra=\big\la\nabla_xG(\zb)^\ast\hat\mu,u+\hat u\big\ra=\skalp{\xba,u+\hat u},
\end{eqnarray*}
and thus $s:=u+\hat u\in\K_{\Gamma(\pb)}(\xb,\xba)$. Next we claim that $\skalp{\nabla_xG(\zb)^\ast\mu,s}=0$. Indeed, supposing the contrary implies that $\skalp{\nabla_xG(\zb)^\ast\mu,s}=\skalp{\mu,\nabla_xG(\zb)s}<0$. Consider now $\bar\mu\in\ri K^\ast$ with $\nabla_x G(\zb)^\ast\bar\mu=\xba$, we get $\bar\mu-t\mu\in K^\ast$ for all $t>0$ sufficiently small and thus arrive at
\begin{eqnarray*}
0\ge\big\la\bar\mu-t\mu,\nabla_xG(\zb)s\big\ra=\big\la\xba,s\big\ra-t\big\la\mu,\nabla_xG(\zb)s\big\ra=-t\big\la\mu,\nabla_xG(\zb)s\big\ra,
\end{eqnarray*}
a contradiction that verifies our claim. It further implies that
\begin{eqnarray*}
0=\big\la\nabla_xG(\zb)^\ast\mu,u+\hat u\big\ra=\big\la\mu,\nabla_xG(\zb)(u+\hat u)\big\ra=\big\la\mu,\nabla_xG(\zb)(u+\hat u)+\th\big\ra=\big\la\mu,\nabla G(\zb)v\big\ra=\big\la\nabla G(\zb)^\ast\mu,v\big\ra.
\end{eqnarray*}
Taking finally into account that $\skalp{\nabla G(\zb)^\ast\hat\mu,v}=0$, we obtain $\skalp{w^\ast,v}=0$ with
\begin{eqnarray*}
w^\ast:=\nabla G(\zb)^\ast(\mu+\alpha\hat\mu)\in\nabla G(\zb)^\ast K^\ast+[\nabla G(\zb)^\ast\hat\mu]\subset\K_{{\rm gph}\,\Gamma}\big(\zb,\nabla G(\zb)^\ast\hat\mu\big)^\ast,
\end{eqnarray*}
which yields $w^\ast\in N_{\K_{{\rm gph}\,\Gamma}(\zb,\nabla G(\zb)^\ast\hat\mu)}(v)$ and thus completes the proof due to $\pi_{X^\ast}(w^\ast)=w_{x^\ast}^\ast$. $\h$\vspace*{0.03in}

Now we are ready to present the aforementioned result of the graphical derivative evaluation without any polyhedrality assumption on the reduced system.\vspace*{-0.05in}

\begin{Theorem}[\bf evaluation of the graphical derivative via the reduction data under strict complementarity]\label{CorStrictComplAux} In addition to the standing assumptions, suppose that the reduced system $G(p,x)\in K$ enjoys the Robinson stability property at $(\pb,\xb)$ and that the strict complementarity condition holds for the system $G(\pb,x)\in K$ at $(\xb,\xba)$. Then for every direction $v=(q,u)\in\Pp\times X$ we have
\begin{eqnarray*}
D\Psi(\pb,\xb,\xba)(v)&\subset&\Big\{\nabla\big(\nabla_x G(\cdot)^\ast\hat\mu\big)(\zb)v+\pi_{X^\ast}\big(N_{\K_{{\rm gph}\,\Gamma}(\zb,z^\ast)}(v)\big)\Big|\\
&&\qquad z^\ast\in N_{T_{{\rm gph}\,\Gamma}(\zb)}(v)\cap\pi_{X^\ast}^{-1}(\xba),\;\hat\mu\in\M(\zb,z^\ast)\Big\}.
\end{eqnarray*}
If furthermore $\nabla^2G(\zb)(v,v)\in\rge\big(\nabla G(\zb)\big)+\lin K$, then this inclusion becomes an equation.
\end{Theorem}\vspace*{-0.05in}
{\bf Proof.} Note first the assumed RS ensures the validity of MSCQ in Lemma~\ref{LemStrictCompl}. Then the claimed inclusion is an immediate consequence of that lemma together with Lemma~\ref{PropTanDir}(i). To verify the equality therein, just note that for all $\mu_1,\mu_2\in\M(\zb,z^\ast)$ we have
\begin{eqnarray*}
\mu_1-\mu_2\in\ker\nabla G(\zb)^\ast\cap(K^\ast-K^\ast)=\big(\rge\nabla G(\zb)+\lin K\big)^\perp\subset[\nabla^2G(\zb)(v,v)]^\perp,
\end{eqnarray*}
which implies that $\big\la\mu_1-\mu_2,\nabla^2G(\zb)(v,v)\big\ra=0$ and $\M(\zb,z^\ast)=\M(\zb,z^\ast;v)$. Hence the asserted equality follows from inclusion  \eqref{RemConv_w_k} in Theorem~\ref{PropSuffTanDir}. $\h$\vspace*{-0.15in}

\section{Computing the Graphical Derivative of the Normal Cone Mapping}\sce\vspace*{-0.1in}

The main intention of this section is to express the computation formulas for our major second-order object $D\Psi$ via the original data of PCS \eqref{EqConstrSystem}. To proceed, observe first that the surjectivity of the derivative $\nabla h(g(\zb))$ in the reduction procedure ensures by the standard chain rule of variational analysis that $N_C(g(\zb))=\nabla h(g(\zb))^\ast K^\ast$ and that the set $\big(\nabla h(g(\zb))^\ast\big)^{-1}(\lambda)$ is a singleton for every $\lambda\in N_C(g(\zb))$. Furthermore, we have $\nabla G(\zb)=\nabla h(g(\zb))\circ\nabla g(\zb)$ and
\begin{eqnarray}
\nonumber\big\la\mu,\nabla^2G(\zb)(v,v)\big\ra&=&\big\la\mu,\nabla h\big(g(\zb)\big)\nabla^2g(\zb)(v,v)+\nabla^2h\big(g(\zb)\big)\big(\nabla g(\zb)v,\nabla g(\zb)v)\big\ra\\
\label{EqChain1}&=&\big\la\nabla h\big(g(\zb)\big)^\ast\mu,\nabla^2g(\zb)(v,v)\big\ra+\nabla^2\big\la\mu,h(\cdot)\big\ra\big(g(\zb)\big)\big(\nabla g(\zb)v,\nabla g(\zb)v\big)\\
\nonumber\nabla\big(\nabla_xG(\cdot)^\ast\mu\big)(\zb)v&=&\nabla\big(\nabla_x g(\cdot)^\ast\nabla h\big(g(\zb)\big)^\ast\mu\big)v+\nabla_xg(\zb)^\ast\big(\nabla h(\cdot)^\ast\mu\big)\big(g(\zb)\big)\nabla g(\zb)v\\
\label{EqChain2}&=&\nabla\big(\nabla_x g(\cdot)^\ast\nabla h\big(g(\zb)\big)^\ast\mu\big)v+\nabla_xg(\zb)^\ast\nabla^2\big\la\mu,h(\cdot)\big\ra\big(g(\zb)\big)\nabla g(\zb)v
\end{eqnarray}
for every $\mu\in E^\ast$ and every $v\in Z$. For each pair $z^\ast\in N_{{\rm gph}\,\Gamma}(\zb)$ we define the {\em multiplier set} by
\begin{eqnarray*}
\Lambda(\zb,z^\ast):=\big\{\lambda\in N_C(g(\zb))\big|\;\nabla g(\zb)^\ast\lambda=z^\ast\big\}
\end{eqnarray*}
and, given any critical direction $v\in\K_{{\rm gph}\,\Gamma}(\zb,z^\ast)$, the {\em directional multiplier set} by
\begin{eqnarray*}
\Lambda(\zb,z^\ast;v):=\argmax_{\lambda\in\Lambda(\zb,\zb^\ast)}\left(\big\la\lambda,\nabla^2g(\zb)(v,v)\big\ra+\nabla^2\big\la\big(\nabla h(g(\zb)\big)^\ast\big)^{-1}(\lambda),h(\cdot)\big\ra\big(g(\zb)\big)\big(\nabla g(\zb)v,\nabla g(\zb)v\big)\right).
\end{eqnarray*}
Then we have $\Lambda(\zb,z^\ast)=\nabla h\big(g(\zb)\big)^\ast\M(\zb,z^\ast)$ and $\Lambda(\zb,z^\ast;v)=\nabla h\big(g(\zb)\big)^\ast\M(\zb,z^\ast;v)$.

The following lemma gives us an expression of the characteristic parameter subspace $\Pp$ from \eqref{char} via the given data of original parametric constraint system \eqref{EqConstrSystem}.\vspace*{-0.07in}

\begin{lemma}[\bf original data description of the characteristic parameter subspace]\label{LemPp} Under the standing assumptions made we have the representation
\begin{equation}\label{EqPp}
\Pp=\big\{q\in P\big|\;\nabla_p g(\zb)q\in\rge\nabla_x g(\zb)+\lin T_C\big(g(\zb)\big)\big\}.
\end{equation}
\end{lemma}\vspace*{-0.05in}
{\bf Proof.} The surjectivity of $\nabla h(g(\zb))$ ensures the equality $T_C(g(\zb))=\{v\mv\nabla h(g(\zb))v\in K\}$. Considering now $q\in\tilde\Pp:=\{q\mv \nabla_pg(\zb)q\in\rge\nabla_x g(\zb)+\lin T_C(g(\zb))\}$, we find $u\in X$ and $v\in\lin T_C(g(\zb))$ with $\nabla _p g(\zb)q=\nabla_x g(\zb)u+v$, which implies the relationships
\begin{eqnarray*}
\nabla_pG(\zb)q=\nabla h\big(g(\zb)\big)\nabla_pg(\zb)q=\nabla h\big(g(\zb)\big)\nabla_xg(\zb)u+\nabla h\big(g(\zb)\big)v=\nabla_xG(\zb)u+\nabla h\big(g(\zb)\big)v.
\end{eqnarray*}
By $v\in\lin T_C(g(\zb))$ we have $\pm v\in T_C(g(\zb)$ and hence $\pm\nabla h(g(\zb)v\in K$, which yields $\nabla h(g(\zb)v\in\lin K$. This gives us $q\in\Pp$ and verifies the inclusion $\tilde\Pp\subset\Pp$. To prove the opposite one, pick $q\in\Pp$, $u\in X$, and $w\in\lin K$ such that $\nabla_pG(\zb)q=\nabla_xG(\zb)u+w$. The surjectivity of $\nabla h(g(\zb))$ allows us to find $v$ with $w=\nabla h(g(\zb))v$, and so $v\in\lin T_C(g(\zb))$ due to $\pm w=\nabla h(g(\zb)(\pm v)\in K$. This tells us that $\nabla h(g(\zb))(\nabla_pg(\zb)q-\nabla_xg(\zb)u-v)=0$ implying the inclusion
\begin{eqnarray*}
\nabla_pg(\zb)q-\nabla_xg(\zb)u-v\in\ker\nabla h\big(g(\zb)\big)\subset\lin T_C\big(g(\zb)\big),
\end{eqnarray*}
which verifies that $\nabla_pg(\zb)q\in\nabla_xg(\zb)u+v+\lin T_C(g(\zb))=\nabla_xg(\zb)u+\lin T_C(g(\zb))$ and so $ q\in\tilde\Pp$. It yields $\Pp\subset\tilde\Pp$ and thus completes the proof of the lemma. $\h$

The next lemma reveals that both MSCQ and RS are invariant with respect to the above reduction.\vspace*{-0.07in}

\begin{lemma}[\bf MSCQ and RC conditions under reduction]\label{LemRS_MSCQ} The following assertions hold:

{\bf(i)} If MSCQ fulfills for the original system $g(z)\in C$ at $\zb$, then it also fulfills for the reduced system $G(z)\in K$ at the same point.

{\bf(ii)} If RS fulfills for the original system $g(p,x)\in C$ at $(\pb,\xb)$, it also fulfills for the reduced system $G(p,x)\in K$ at the same point.
\end{lemma}\vspace*{-0.07in}
{\bf Proof.} We verify assertions (i) and (ii) in a parallel way. Consider neighborhoods $Q$ of $\pb$ and $U$ of $\xb$ together with a positive number $\kappa$ such that the estimate ${\rm dist}(z;g^{-1}(C))\le\kappa\,{\rm dist}(g(z),C)$ holds for all $z\in Q\times U$ in case (i) and the estimate
${\rm dist}(x;\Gamma(p))\le\kappa\,{\rm dist}(g(p,x),C)$ holds for all $(p,x)\in Q\times P$ in case (ii). The surjectivity of $\nabla h(\yb)$ is equivalent to the metric regularity of the mapping $h(\cdot)-K$ around $(g(\zb),0)$ (see, e.g., \cite[Theorem~4.18]{M06}) and thus we can find $\rho,\beta>0$ such that
\begin{eqnarray*}
{\rm dist}\big(y;h^{-1}(K)\big)\le\beta\,{\rm dist}\big(h(y);K\big)\;\mbox{ for all }\;y\in\B_\rho\big(g(\zb)\big).
\end{eqnarray*}
Choose $\rho$ to be so small that $C\cap\B_\rho(g(\zb))=h^{-1}(K)\cap\B_\rho(g(\zb))$ and then find $r>0$ with $\B_r(\zb)\subset g^{-1}(\B_\rho(g(\zb)))$.
Set $\tilde\rho:=\frac 12\min\{\frac r\kappa,\rho\}$ and consider any vector
\begin{eqnarray*}
z=(p,x)\in(Q\times U)\cap g^{-1}\big(\B_{\tilde\rho}(g(\zb))\big)\cap\B_{r/2}(\zb)
\end{eqnarray*}
from the neighborhood of $\zb$. Since ${\rm dist}\big(g(z);C\big)\le\norm{g(z)-g(\zb)}\le\tilde\rho$, there is $c\in C$ with ${\rm dist}(g(z);C)=\norm{g(z)-c}$ and $\norm{c-g(\zb)}\le\norm{c-g(z)}+\norm{g(z)-g(\zb)}\le\frac\rho 2+\frac\rho 2=\rho$. Thus
\begin{eqnarray*}
{\rm dist}(g(z);C)={\rm dist}\big(g(z);C\cap\B_\rho(g(\zb))\big)={\rm dist}\big(g(z);h^{-1}(K)\cap\B_\rho(g(\zb))\big)={\rm dist}\big(g(z);h^{-1}(K)\big).
\end{eqnarray*}
Dealing now with the MSCQ case (i), we find $\Tilde z\in g^{-1}(C)-z$ with $\norm{\tilde z}={\rm dist}(z;g^{-1}(C))\le\kappa\,{\rm dist}(g(z);C)\le\kappa\, \tilde\rho\le\frac r2$. It shows that $\norm{z+\tilde z-\zb}\le\norm{z-\zb}+\norm{\tilde z}\le\frac r2+\frac r2=r$ and hence $g(z+\tilde z)\in C\cap \B_\rho(g(\zb))=h^{-1}(K)\cap\B_\rho(g(\zb))$. This implies that $G(z+\tilde z))=h(g(z+\tilde z))\in K$ and
\begin{eqnarray*}
{\rm dist}\big(z;G^{-1}(K)\big)\le\norm{\tilde z}={\rm dist}\big(z;g^{-1}(C)\big)\le\kappa{\rm dist}\big(g(z);C\big)=\kappa{\rm dist}\big(g(z);h^{-1}(K)\big)\le\kappa\beta{\rm dist}\big(h(g(z));K\big),
\end{eqnarray*}
which verifies the metric subregularity of $G(\cdot)-K$ at $(\zb,0)$.

In the RS case (ii) we find $\tilde x\in\Gamma(p)-x$ with  $\norm{\tilde x}={\rm dist}(x;\Gamma(p))\le\kappa{\rm dist}(g(z);C)\le\kappa\tilde\rho\le\frac r2$. Setting $\tilde z:=(0,\tilde x)$ gives us $\norm{(z+\tilde z-\zb}\le\norm{z-\zb}+\norm{\tilde z}\le\frac r2+\frac r2=r$, and so $g(z+\tilde z)\in C\cap \B_\rho(g(\zb))=h^{-1}(K)\cap\B_\rho(g(\zb))$. Hence we get $G(z+\tilde z))=h(g(z+\tilde z))\in K$ and
\begin{eqnarray*}
{\rm dist}\big(x;G(\cdot,p)^{-1}(K)\big)\le\norm{\tilde x}={\rm dist}\big(x;\Gamma(p)\big)\le\kappa{\rm dist}\big(g(z);C\big)=\kappa{\rm dist}\big(g(z);h^{-1}(K)\big)\le\kappa\beta{\rm dist}\big(h(g(z));K\big),
\end{eqnarray*}
and verifies the Robinson stability of the reduced system $G(p,x)\in K$ at $(\pb,\xb)$. $\h$\vspace*{0.06in}

To formulate the first major result of this section on computing the graphical derivative of $\Psi$ in \eqref{EqNormalConeMap}, for every $\lambda\in N_C(g(\zb))$
define the linear mapping ${\cal H}_\lambda\colon Z\to X^\ast$ by
\begin{eqnarray}\label{H1}
{\cal H}_\lambda v:=\nabla_xg(\zb)^\ast\nabla^2\big\la\big(\nabla h(g(\zb))^\ast\big)^{-1}\lambda,h(\cdot)\big\ra\big(g(\zb)\big)\nabla g(\zb)v.
\end{eqnarray}
Observe that, although mapping \eqref{H1} (as well as its modifications considered below) is mainly constructed via the original data of the constraint system \eqref{EqConstrSystem}, it still includes the $C^2$-smooth mapping $h\colon V\to E$ that reduces the given constraint set $C$ to the cone $K$ in \eqref{DefConeReducible}, which is not in the picture anymore. However, the mapping $h$ can be directly expressed entirely via the given data of \eqref{EqConstrSystem} in the most interesting settings. In particular, $h$ is the identity mapping if $C$ is a closed convex cone itself and $g(\op,\ox)=0\in C$. In the case of the ice-cream/Lorentz cone $C\subset Y=\R^l$ we have that $h(y)$ is the identity mapping for $y=0$ and $h(y)=y_1^2+\ldots+y^2_{l-1}-y^2_l$ otherwise, while for the ice-cream cone products the mapping $h(y)$ at nonzero points is composed by the quadratic forms depending on the values of the components at the reference points. We refer the reader to \cite[Example~3.140]{BonSh00} for the explicit calculation of the reduction mapping $h$ in the important case of the SDP cone. Some other examples of calculating reduction mappings can also be found in \cite{BonSh00}.

Note furthermore that geometrically the mapping ${\cal H}_\lm$ in \eqref{H1} describes the {\em curvature} of the set $C$ at the point $g(\op,\ox)$. In particular, there is no curvature if $g(\op,\ox)=0$ and so any tangent to the convex cone $C$ at the origin is contained in $C$. But it is not the case if $g(\op,\ox)\ne 0$ and the set $C$ is not polyhedral. The following theorem reflects all of this in the second-order computation formulas.\vspace*{-0.05in}

\begin{Theorem}[\bf computing the graphical derivative via the given data of PCS]\label{CorMainTh}
In addition to the standing assumptions, suppose that the original parametric constraint system $g(p,x)\in C$ enjoys the RS property at $(\pb,\xb)$ and that either $\Pp=P$ or the cone $K$ is polyhedral. Then we have
\begin{eqnarray*}
T_{{\rm gph}\,\Psi}(\pb,\xb,\xba)&=&\Big\{(q,u,u^\ast)\Big|\,\exists\,z^\ast\in N_{T_{{\rm gph}\,\Gamma}(\zb)}(q,u)\cap \pi_{X^\ast}^{-1}(\xba),\;\hat\lambda\in\Lambda\big(\zb,z^\ast;(q,u)\big)\\
&&\quad\mbox{s.t. }\;u^\ast\in\Big(\nabla\big(\nabla_x g(\cdot)^\ast\hat\lambda\big)(\zb)+{\cal H}_{\hat\lambda}\Big)(q,u)+\pi_{X^\ast}\big(N_{\K_{{\rm gph}\,\Gamma}(\zb,z^\ast)}(q,u)\big)\Big\}.
\end{eqnarray*}
Consequently, for all $v=(q,u)\in Z$ the graphical derivative of $\Psi$ is computed by
\begin{eqnarray*}
D\Psi(\pb,\xb,\xba)(v)&=&\Big\{\Big(\nabla\big(\nabla_x g(\cdot)^\ast\hat\lambda\big)(\zb)+{\cal H}_{\hat\lambda}\Big)v+\pi_{X^\ast}\big(N_{\K_{{\rm gph}\,\Gamma}(\zb,z^\ast)}(v)\big)\Big|\\
&&\quad z^\ast\in N_{T_{{\rm gph}\,\Gamma}(\zb)}(v)\cap\pi_{X^\ast}^{-1}(\xba),\;\hat\lambda\in\Lambda(\zb,z^\ast;v)\Big\}.
\end{eqnarray*}
\end{Theorem}\vspace*{0.05in}
{\bf Proof.} We deduce this statement from Theorem \ref{ThMainTh} and Lemma~\ref{LemRS_MSCQ} by putting $\nabla h(g(\zb))^\ast\mu=\lambda$ therein and using the above formulas \eqref{EqChain1} and \eqref{EqChain2}.$\h$\vspace*{0.05in}

The assumption $\Pp=P$ is obviously fulfilled when $\nabla_pg(\pb,\xb))=0$ (``weak parameterization") that holds, in particular, in the nonparametric setting of $g(p,\cdot)=g(\pb,\cdot)$ for all $p$. This leads us to the following striking consequence of Theorem~\ref{CorMainTh}, which significantly extends  all the known results in this direction obtained for particular classes of constraint systems.\vspace*{-0.05in}

\begin{corollary}[\bf graphical derivative computation for nonparametric constraint systems]\label{nonpar} Consider the constraint system $g(x)\in C$ generated by a $C^2$-smooth mapping $g\colon X\to Y$ and a set $C\subset Y$ that is $C^2$-cone reducible at $g(\xb)$ with $\xb\in\tilde\Gamma:=\{x\in X\mv g(x)\in C\}$. Take $\xba\in N_{\tilde\Gamma}(\xb)$ and assume that MSCQ holds for the system $g(x)\in C$ at $\xb$. Then we have the tangent cone formula
\begin{eqnarray*}
T_{{\rm gph}\,N_{\tilde\Gamma}}(\xb,\xba)=\Big\{(u,u^\ast)\Big|\;\exists\,\hat\lambda\in\Lambda(\xb,\xba;u)\;\mbox{ with }\;
u^\ast\in\Big(\nabla^2\big\la\hat\lambda,g\big\ra(\xb)+\tilde{\cal H}_{\hat\lambda}\Big)u+N_{\K_{\Gamma}(\xb,\xba)}(u)\Big\},
\end{eqnarray*}
where the linear mapping $\tilde{\cal H}_{\hat\lambda}\colon X\to X^\ast$ is defined by
\begin{eqnarray*}
\tilde{\cal H}_{\hat\lambda}u:=\nabla g(\xb)^\ast\nabla^2\big\la(\nabla h(g(\xb))^\ast\big)^{-1}\hat\lambda,h(\cdot)\big\ra\big(g(\xb)\big)\nabla g(\xb)u.
\end{eqnarray*}
Consequently, for all $u\in X$ the graphical derivative of the normal cone mapping is computed by
\begin{eqnarray*}
D N_{\tilde\Gamma}(\xb,\xba)(u)&=&\Big\{\Big(\nabla^2\big\la\hat\lambda,g\big\ra(\xb)+\Tilde{\cal H}_{\hat\lambda}\Big)u+N_{\K_{\tilde\Gamma}(\xb,\xba)}(u)\Big| \hat\lambda\in\Lambda(\xb,\xba;u)\Big\}.
\end{eqnarray*}
\end{corollary}\vspace*{-0.05in}
{\bf Proof.} This is a clear consequence of Theorem~\ref{CorMainTh}, where RC reduces to MSCQ and the mapping ${\cal H}_\lm$ from \eqref{H1} agrees with $\tilde{\cal H}_\lm$ in the nonparametric setting under consideration. $\h$\vspace*{0.05in}

The results obtained in Corollary~\ref{nonpar} extend those derived in \cite{CH17,GO16} for polyhedral systems (when the terms with $\tilde{\cal H}$ disappear in the above formulas) and in \cite{HMS17} for the case of the (single) ice-cream cone $C$, where the curvature term $\tilde{\cal H}$ reduces to the one in \cite[Theorem~5.1]{HMS17}, which was first introduced in \cite{BonRa05} for different purposes in second-order cone/ice-cream programming.\vspace*{0.02in}

The last consequence of Theorem~\ref{CorMainTh} in this section concerns the general parametric setting of \eqref{EqConstrSystem} under the assumptions of the theorem and reveals relationships between the graphical derivatives of $\Psi$ and of the normal cone mappings for the graph of $\Gamma$ computed via the data of \eqref{EqConstrSystem}.\vspace*{-0.05in}

\begin{corollary}[\bf relationships for graphical derivatives of parametric constraint systems with $\Pp=P$]\label{gd-relat} In the setting and under the
assumptions of Theorem~{\rm\ref{CorMainTh}} we have the formula
\begin{eqnarray*}
T_{{\rm gph}\,\Psi}(\pb,\xb,\xba)&=&\bigcup_{z^\ast\in N_{{\rm gph}\,\Gamma}(\zb)\cap\pi_{X^\ast}^{-1}(\xba)}\Big\{\big(v,\pi_{X^\ast}(v^\ast)\big)\Big|\;(v,v^\ast)\in T_{{\rm gph}\,N_{{\rm gph}\,\Gamma}}(\zb,z^\ast)\Big\}
\end{eqnarray*}
and correspondingly the graphical derivative relationship valid for every $v=(q,u)$:
\begin{eqnarray*}
D\Psi(\pb,\xb,\xba)(v)=\bigcup_{z^\ast\in N_{T_{{\rm gph}\,\Gamma}(\zb)}(v)\cap\pi_{X^\ast}^{-1}(\xba)}\pi_{X^\ast}\big(DN_{{\rm gph}\,\Gamma}(\zb,z^\ast)(v)\big),
\end{eqnarray*}
where the tangent cone and the graphical derivative on the right-hand sides above are computed by
\begin{eqnarray*}
T_{{\rm gph}\,N_{{\rm gph}\,\Gamma}}(\zb,z^\ast)=\Big\{(v,v^\ast)\Big|\;\exists\,\hat\lambda\in\Lambda(\zb,z^\ast;v)\;\mbox{ s.t. }\;
v^\ast\in\Big(\nabla^2\big\la\hat\lambda,g\big\ra(\zb)+\tilde{\cal H}_{\hat\lambda}\Big)(v)+N_{\K_{{\rm gph}\,\Gamma}(\zb,z^\ast)}(v)\Big\},
\end{eqnarray*}
\begin{eqnarray*}
DN_{{\rm gph}\,\Gamma}(\zb,z^\ast)(v)&=&\Big\{\Big(\nabla^2\big\la\hat\lambda,g\big\ra(\zb)+\tilde{\cal H}_{\hat\lambda}\Big)(v)+N_{\K_{{\rm gph}\,\Gamma}(\zb,z^\ast)}(v)\Big|\,\hat\lambda\in\Lambda(\zb,z^\ast;v)\Big\}
\end{eqnarray*}
with $\tilde{\cal H}_{\hat\lambda}v:=\nabla g(\zb)^\ast\nabla^2\big\la\big(\nabla h(g(\zb))^\ast\big)^{-1}\hat\lambda,h(\cdot)\big\ra\big(g(\zb)\big)\nabla g(\zb)v$.
\end{corollary}\vspace*{-0.05in}
{\bf Proof.} It follows from the definitions and the results of the theorem by applying the corresponding formulas of Corollary~\ref{nonpar} to the set $\tilde\Gamma=\gph\Gamma$ therein. $\h$\vspace*{0.05in}

Note that in the case where $g(p,x)$ {\em strongly} depends on the parameter (``strong parameterization") meaning that the derivative $\nabla_p g(\op,\ox)$ is surjective, the condition $\Pp=P$ amounts to the requirement
\begin{eqnarray*}
\rge\nabla_xg(\oz)+\lin T_C\big(g((\oz)\big)=Y,
\end{eqnarray*}
which is the well-known {\em nondegeneracy condition} for PCS \eqref{EqConstrSystem} at $(\op,\ox)$; see \cite{BonSh00}. The next major result is a version of the reduction Theorem~\ref{CorStrictComplAux} in terns of the original PCS data without imposing either the assumption on $\Pp=P$ or on the polyhedrality of $K$ while replacing them by {\em strict complementarity}. However, the precise graphical derivative computation (versus an upper estimate useful for its own sake) is obtained under an additional condition, which surely holds for the case of strong parameterization while being far enough from nondegeneracy even in this case.\vspace*{-0.05in}

\begin{Theorem}[\bf graphical derivative for the original PCS under strict complementarity]\label{CorStrictCompl} Together with the standing assumptions on the original PCS \eqref{EqConstrSystem} as well as the Robinson stability imposed on it at $(\pb,\xb)$, suppose that the strict complementarity condition holds for the system $g(\pb,x)\in C$ at $(\xb,\xba)$. Then for every direction $v=(q,u)\in\Pp\times X$ with parameters from \eqref{EqPp} we have the inclusion
\begin{eqnarray*}
D\Psi(\pb,\xb,\xba)(v)&\subset&\Big\{\Big(\nabla\big(\nabla_x g(\cdot)^\ast\hat\lambda\big)(\zb)+{\cal H}_{\hat\lambda}\Big)v+\pi_{X^\ast}\big(N_{\K_{{\rm gph}\,\Gamma}(\zb,z^\ast)}(v)\big)\Big|\\
&&\qquad z^\ast\in N_{T_{{\rm gph}\,\Gamma}(\zb)}(v)\cap\pi_{X^\ast}^{-1}(\xba),\;\hat\lambda\in\Lambda(\zb,z^\ast)\Big\}
\end{eqnarray*}
with ${\cal H}_{\lambda}$ taken from \eqref{H1}. This inclusion holds as equality if in addition we assume that
\begin{equation}\label{EqConstLambda}
\nabla h\big(g(\zb)\big)\nabla^2g(\zb)(v,v)+\nabla^2h\big(g(\zb)\big)\big(\nabla g(\zb)v,\nabla g(\zb)v)\in\rge\big(\nabla h\big(g(\zb)\big)\circ\nabla g(\zb)\big)+\lin K,
\end{equation}
which is satisfied, in particular, when the derivative $\nabla g(\zb)$ is surjective.
\end{Theorem}\vspace*{-0.05in}
{\bf Proof.} Under the imposed strict complementarity condition of the theorem we take $\bar\lambda\in\ri N_C(g(\zb))$ and deduce from basic convex analysis that
\begin{eqnarray*}
\bar\mu:=\big(\nabla h(g(\zb))^\ast\big)^{-1}\bar\lambda\in\ri\Big(\big(\nabla h(g(\zb))^\ast\big)^{-1}N_C(g(\zb))\Big)=\ri K^\ast.
\end{eqnarray*}
This shows that the strict complementarity condition also holds for the reduced system $G(\pb,x)\in K$ at $(\xb,\xba)$. Thus it remains to apply Theorem~\ref{CorStrictComplAux} with the usage of the relationships between the original and reduced constraint systems that are discussed above. $\h$\vspace*{0.03in}

The concluding result of this section does not impose either $\Pp=P$, or polyhedrality, or strict complementarity assumptions on PCS while establishing an upper estimate for a modified graphical derivative of the normal cone mapping $\Psi$. It is useful for stability applications in the next section. The {\em modified graphical derivative} $\tilde D\Psi(\pb,\xb,\xba)\colon P\times X\tto X^\ast$ of the mapping $\Psi$ from \eqref{EqNormalConeMap} is defined by
\begin{align}\label{mgd}
\tilde D\Psi(\pb,\xb,\xba)(q,u):=\Big\{u^\ast\Big|&\exists\,t_k\downarrow 0,\;(q_k,u_k,u_k^\ast)\to(q,u,u^\ast)\;\mbox{ with }\\\nonumber
&\xba+t_ku_k^\ast\in\Psi(\pb+t_kq_k,\xb+t_ku_k),\;\limsup_{k\to\infty}\frac{\dist(q_k-q;\Pp)}{t_k}<\infty\Big\}.
\end{align}\vspace*{-0.25in}

\begin{proposition}[\bf evaluation of the modified graphical derivative]\label{CorRemCases} Add only the validity of RS for $g(p,x)\in C$ at $(\op,\ox)$ to our standing assumptions. Then for all $v=(q,u)\in P\times X$ we have the inclusions
\begin{eqnarray*}
\begin{array}{ll}
\tilde D\Psi(\pb,\xb,\xba)(q,u)\subset\Big\{\Big(\nabla\big(\nabla_x g(\cdot)^\ast\hat\lambda\big)(\zb)+{\cal H}_{\hat\lambda}\Big)v+w_{x^\ast}^\ast\Big|\;
z^\ast\in N_{T_{{\rm gph}\,\Gamma}(\zb)}(v)\cap\pi_{X^\ast}^{-1}(\xba),\;\hat\lambda\in\Lambda(\zb,z^\ast;v),\\
\exists w_k^\ast\in\K_{{\rm gph}\,\Gamma}(\zb,z^\ast)\big)^\ast\;\mbox{ s.t. }\;
\disp w_{x^\ast}^\ast=\lim_{k\to\infty}\pi_{X^\ast}(w_k^\ast),\;\lim_{k\to\infty}\skalp{w_k^\ast,v}=0\Big\}\subset D\Psi(\pb,\xb,\xba)(q,u)\\
\subset\Big\{\Big(\nabla\big(\nabla_x g(\cdot)^\ast\hat\lambda\big)(\zb)+{\cal H}_{\hat\lambda}\Big)v+\K_{\Gamma(\pb)}(\xb,\xba)^\ast\Big|\;\hat\lambda\in N_{T_C(g(\zb))}\big(\nabla g(\zb)v\big),\nabla_xg(\zb)^\ast\hat\lambda=\xba\Big\}.
\end{array}
\end{eqnarray*}
\end{proposition}\vspace*{-0.05in}
{\bf Proof.} It follows from Lemma~\ref{PropTanDir}(iii) and the relationships between the original and reduced parametric constraint systems revealed above. $\h$\vspace*{-0.15in}

\section{Applications to Parametric Variational Systems}\sce\vspace*{-0.1in}

In this section we study stability properties of {\em parametric variational systems} (PVS) of the type
\begin{equation}\label{EqSolMap}
S(p):=\big\{x\in X\big|\;0\in f(p,x)+N_{\Gamma(p)}(x)\big\},
\end{equation}
which are given as solution maps to the parameter-dependent generalized equations $0\in f(p,x)+N_{\Gamma(p)}(x)$. The set-valued mapping $\Gamma\colon P\tto X$ in \eqref{EqSolMap} is taken from \eqref{EqGamma} under the standing assumptions imposed on it in Section~2, while the single-valued mapping $f\colon P\times X\to X^\ast$ is assumed to be continuously differentiable around the reference point. We are not going to mention these assumptions in the rest of the section. Observe that in the case where $f(p,x)=\nabla_x\ph(p,x)$ for some differentiable function $\ph\colon P\times X\to\R$, the variational system $S(p)$ in \eqref{EqSolMap} describes the collections of points satisfying the basic first-order necessary optimality conditions in the parametric optimization problem:
\begin{equation*}
\min_x\ph(p,x)\quad\mbox{subject to}\quad g(p,x)\in C.
\end{equation*}

Let us fix the reference point $\zb=(\pb,\xb)\in\gph S$ and we consider the associated set of multipliers
\begin{eqnarray}\label{mult1}
\Lb:=\big\{\lambda\in N_C(g(\zb))\big|\;f(\zb)+\nabla_xg(\zb)^\ast\lambda=0\big\}.
\end{eqnarray}
Given $\lambda\in N_C(g(\zb))$, define the linear operator ${\cal H}_\lambda^x\colon X\to X^\ast$ by
\begin{eqnarray}\label{H2}
{\cal H}_\lambda^x u:=\nabla_xg(\zb)^\ast\nabla^2\big\la\big(\nabla h(g(\zb))^\ast\big)^{-1}\lambda,h(\cdot)\big\ra\big(g(\zb)\big)\nabla_xg(\zb)u,
\end{eqnarray}
which is a partial version of \eqref{H1} satisfying the condition ${\cal H}_\lambda^x u={\cal H}_\lambda(0,u)$, and then define yet another linear operator $\F_\lambda\colon X\to X^\ast$ by
\begin{eqnarray}\label{F}
\F_\lambda u:=\nabla_x f(\pb,\xb)u+\nabla_x^2\big\la\lambda,g(\cdot)\big\ra(\zb)u+{\cal H}_\lambda^x u.
\end{eqnarray}

Our first result here concerns not isolated calmness (which is Lipschitzian in nature; see below) but its {\em H\"olderian} version that is of its own interest while being (together with its proof) of a crucial importance in verifying isolated calmness in what follows.\vspace*{-0.07in}

\begin{Theorem}[\bf H\"olderian isolated calmness of PVS]\label{ThHoelder} Assume that PCS \eqref{EqConstrSystem} enjoys the Robinson stability property at $(\pb,\xb)$ and that the following condition is satisfied:\vspace*{-0.1in}
\begin{equation}\label{EqSuffCond}
\parbox{\myparboxwidth}{For every $0\not=u\in\K_{\Gamma(\pb)}\big(\xb,-f(\zb)\big)$, every $\lambda\in\Lb$, and every $\hat\lambda\in\Lambda\big(\zb,\nabla g(\zb)^\ast\lambda;(0,u)\big)$ we have either $\skalp{\F_{\hat\lambda}u,u}>0$ or $\skalp{\F_{\hat\lambda}u,\tilde u}<0$ with some $\tilde u\in\K_{\Gamma(\pb)}\big(\xb,-f(\zb)\big)$.}\vspace*{-0.1in}
\end{equation}
Then there are a constant $\ell>0$ and neighborhoods $U$ of $\xb$ and $Q$ of $\pb$ such that
\begin{eqnarray}\label{hol}
S(p)\cap U\subset\big\{\xb\big\}+\ell\big(\norm{p-\pb}+\sqrt{{\rm dist}(p-\pb;\Pp)}\big)\B_X\;\mbox{ whenever }\;p\in Q.
\end{eqnarray}
\end{Theorem}\vspace*{-0.1in}
{\bf Proof.} To verify \eqref{hol}, suppose on the contrary that there is a sequence $(p_k,x_k)\longsetto{{{\rm gph}\,S}}(\pb,\xb)$ with
\begin{eqnarray*}
\norm{x_k-\xb}>k\big(\norm{p_k-\pb}+\sqrt{{\rm dist}(p_k-\pb;\Pp)}\big)\;\mbox{ for all }\;k\in\N.
\end{eqnarray*}
Define $t_k:=\norm{x_k-\xb}$, $q_k:=(p_k-\pb)/t_k$, and $u_k:=(x_k-\xb)/t_k$ and then get $\norm{q_k}<\frac 1k$ with $u_k\to u\not=0$ along a subsequence. Furthermore, it follows that
\begin{eqnarray*}
{\rm dist}(q_k-0;\Pp)=\frac{{\rm dist}(p_k-\pb;\Pp)}{t_k}<\frac{(t_k/k)^2}{t_k}=\frac{t_k}{k^2}\;\mbox{ for all }\;k\in\N,
\end{eqnarray*}
and hence $\lim_{k\to\infty}{\rm dist}(q_k-0;\Pp)/t_k=0$ as $k\to\infty$. By $x_k\in S(p_k)$ it tells us that
\begin{eqnarray*}
f(\pb+t_kq_k,\xb+t_ku_k)=-\big(f(\zb)+t_k\nabla_xf(\zb)u_k+\oo(t_k)\big)\in\Psi(\pb+t_k q_k,\xb+t_ku_k),
\end{eqnarray*}
which being combined with $\xba:=-f(\zb)\in\Psi(\pb,\xb)$ brings us to the inclusion
\begin{eqnarray*}
0\in\nabla_xf(\zb)u+\tilde D\Psi\big(\pb,\xb,-f(\zb)\big)(0,u),
\end{eqnarray*}
where the modified graphical derivative $\tilde D$ is taken from \eqref{mgd}. Applying now Proposition~\ref{CorRemCases}, we find $z^\ast\in N_{T_{{\rm gph}\,\Gamma}(\zb)}(0,u)\cap \pi_{X^\ast}^{-1}(\xba)$, $\hat\lambda\in\Lambda(\zb,z^\ast;(0,u))$, and a sequence of
$w_k^\ast\in\big(\K_{{\rm gph}\,\Gamma}(\zb,z^\ast)\big)^\ast$ such that $\lim_{k\to\infty}\skalp{w_k^\ast,(0,u)}=\lim_{k\to\infty}\skalp{\pi_{X^\ast}(w_k^\ast),u}=0$ and
\begin{eqnarray*}
0=\nabla_x f(\zb)u+\Big(\nabla\big(\nabla_x g(\cdot)^\ast \hat\lambda\big)(\zb)+{\cal H}_{\hat\lambda}\Big)(0,u)+\lim_{k\to\infty}\pi_{X^\ast}(w_k^\ast)=\F_{\hat\lambda}u +\lim_{k\to\infty}\pi_{X^\ast}(w_k^\ast)
\end{eqnarray*}
with ${\cal H}_\lambda$ taken from \eqref{H1} and ${\cal F}_\lm$ defined in \eqref{F}. Thus $(0,u)\in T_{{\rm gph}\,\Gamma}(\zb)$, which implies that $\nabla g(\zb)(0,u)=\nabla_xg(\zb)u\in T_C(g(\zb))$ and $u\in T_{\Gamma(\pb)}(\xb)$. Since we clearly have $z^\ast\in N_{T_{{\rm gph}\,\Gamma}(\zb)}(0,u)=N_{{\rm gph}\,\Gamma}(\zb)\cap[(0,u)]^\perp$, there is $\lambda\in N_C(g(\zb))$ satisfying $z^\ast=\nabla g(\zb)^\ast\lambda$. Combining it with $\pi_{X^\ast}(z^\ast)=\nabla_xg(\zb)^\ast\lambda=\xba$ yields $\lambda\in\Lb$. It follows further that $\skalp{z^\ast,(0,u)}=\skalp{\pi_{X^\ast}(z^\ast),u}=\skalp{\xba,u}=0$, and therefore $u\in\K_{\Gamma(\pb)}(\xb,\xba)$. Proposition~\ref{LemNormalCone} allows us to find $\lambda_k\in N_C(g(\zb))$ and $\alpha_k\in\R$ such that
\begin{eqnarray*}\norm{w_k^\ast-(\nabla g(\zb)^\ast\lambda_k+\alpha_k z^\ast)}\le\frac 1k\;\mbox{ for all }\;k\in\N\;\mbox{ and}
\end{eqnarray*}\vspace*{-0.3in}
\begin{eqnarray*}
\disp\lim_{k\to\infty}\pi_{X^\ast}(w_k^\ast)=\lim_{k\to\infty}\pi_{X^\ast}\big(\nabla g(\zb)^\ast\lambda_k+\alpha_k z^\ast)\big)=\lim_{k\to\infty}\big(\nabla_xg(\zb)^\ast\lambda_k+\alpha_k\xba)\big).
\end{eqnarray*}
Using Proposition~\ref{LemNormalCone} again gives us $\lim_{k\to\infty}(\nabla_xg(\zb)^\ast\lambda_k+\alpha_k \xba))\in\K_{\Gamma(\pb)}(\xb,\xba)^\ast$, which implies together with the limiting condition $\lim_{k\to\infty}\skalp{\pi_{X^\ast}(w_k^\ast),u}=0$ that
\begin{eqnarray*}
\disp\lim_{k\to\infty}\pi_{X^\ast}(w_k^\ast)\in\K_{\Gamma(\pb)}(\xb,\xba)^\ast\cap[u]^\perp=\big(\K_{\Gamma(\pb)}(\xb,\xba)+[u]\big)^\ast.
\end{eqnarray*}
Hence $-\F_{\hat\lambda}u\in\big(\K_{\Gamma(\pb)}(\xb,\xba)+[u]\big)^\ast$, which is equivalent to
\begin{eqnarray*}
\big\la-\F_{\hat\lambda}u,\tilde u+\alpha u\big\ra\le 0\;\mbox{ for all }\;\tilde u\in\K_{\Gamma(\pb)}(\xb,\xba)\;\mbox{ and }\;\alpha\in\R.
\end{eqnarray*}
The latter amounts to saying that $\skalp{\F_{\hat\lambda}u,\tilde u}\ge 0$ for all $\tilde u\in\K_{\Gamma(\pb)}$ and that $\skalp{\F_{\hat\lambda}u, u}=0$. This clearly contradicts the assumptions in \eqref{EqSuffCond} and thus completes the proof of the theorem. $\h$\vspace*{0.03in}

To proceed further, for an arbitrary number $\gamma>0$ we define the subset of parameters
\begin{eqnarray*}
\Pp_\gamma:=\big\{p\in P\big|\;{\rm dist}(p-\pb;\Pp)\le\gamma\norm{p-\pb}^2\big\}
\end{eqnarray*}
and deduce from Theorem~\ref{ThHoelder} the inclusion
\begin{eqnarray}\label{hol1}
S(p)\cap U\subset\big\{\xb\big\}+(\ell+\gamma)\norm{p-\pb}\B_X\;\mbox{ for all }\;p\in Q\cap\Pp_\gamma.
\end{eqnarray}
Recall that the mapping $S$ is said to be {\em isolatedly calm} at $(\pb,\xb)$ if there exist
a constant $\ell>0$ and neighborhoods $Q$ of $\pb$ and $U$ of $\xb$ such that
\begin{eqnarray}\label{calm}
S(p)\cap U\subset\big\{\xb\big\}+\ell\norm{p-\pb}\B_X\;\mbox{ for all }\;p\in Q.
\end{eqnarray}
This (Lipschitzian) stability property has been recognized in variational analysis and its applications while being equivalent to the {\em strong metric subregularity} of the inverse $S^{-1}$; see \cite{DoRo14} for more details and references. The result in \eqref{hol1} of Theorem~\ref{ThHoelder} tells us therefore that the {\em restriction} of the solution map \eqref{EqSolMap} to the parameter subset $\Pp_\gamma$ enjoys the isolated calmness property at $(\op,\ox)$.

As an immediate consequence of \eqref{hol}, we get the isolated calmness of $S$ at $(\pb,\xb)$ when $\Pp=P$ under the assumptions of Theorem~\ref{ThHoelder}. Besides it, the next theorem provides other fairly general conditions ensuring the validity of isolated calmness for PVS \eqref{EqSolMap}. Its proof is based on applying the second-order computations obtained above and the {\em graphical derivative criterion} for the isolated calmness property of any closed-graph mapping $S\colon\R^m\tto\R^n$ at $(\op,\ox)\in\gph S$ that reads as
\begin{eqnarray}\label{gd-cr}
DS(\op,\ox)(0)=\{0\}.
\end{eqnarray}
This criterion was explicitly established by Levy \cite{L96} while its derivation can actually be found in Rockafellar \cite{Roc89}; see also \cite[Theorem~4E.1 and Corollary~4E.2]{DoRo14} for more details and discussions.\vspace*{-0.07in}

\begin{Theorem}[\bf sufficient conditions for isolated calmness in PVS]\label{ThCalm} Under the Robinson stability of $g(p,x)\in C$ at $(\op,\ox)$, suppose that the assumptions in one of the following statements are satisfied:

{\bf (i)} In addition to \eqref{EqSuffCond}, either $\Pp=P$ or $K$ is polyhedral.

{\bf (ii)} Strict complementarity holds for the system $g(\pb,x)\in C$ at $(\xb,-f(\zb))$ and\vspace*{-0.1in}
\begin{equation}\label{EqStrongSuffCond}
\parbox{\myparboxwidth}{for all $0\not=u\in\K_{\Gamma(\pb)}\big(\xb,-f(\zb)\big)$ and $\lambda\in\Lb$ from \eqref{mult1}
we have either $\skalp{\F_{\lambda}u,u}>0$ or $\skalp{\F_{\lambda}u,\tilde u}<0$ with some $\tilde u\in\K_{\Gamma(\pb)}\big(\xb,-f(\zb)\big)$.}\vspace*{-0.1in}
\end{equation}
Then PVS \eqref{EqSolMap} enjoys the isolated calmness property \eqref{calm} at $(\op,\ox)$.
\end{Theorem}\vspace*{-0.07in}
{\bf Proof.} First we verify the claimed isolated calmness under the assumptions in (i). Arguing by contradiction and using the graphical derivative criterion \eqref{gd-cr} together with the graphical derivative construction in \eqref{gd}, \eqref{tan} and the form of $S$ in \eqref{EqSolMap}, we find $u\ne 0$ such that
\begin{eqnarray}\label{sum}
0\in\nabla_x f(\zb)u+D\Psi(\zb,\xba)(0,u)\;\mbox{ with }\;\xba:=-f(\zb).
\end{eqnarray}
Then the graphical derivative computation of Theorem~\ref{CorMainTh} gives us vectors $z^\ast\in N_{T_{{\rm gph}\,\Gamma}(\zb)}(0,u)\cap\pi_{X^\ast}^{-1}(\xba)$, $\hat\lambda\in\Lambda(\zb,z^\ast;(0,u))$, and $w^\ast\in N_{\K_{{\rm gph}\,\Gamma}(\zb,z^\ast)}(0,u)$ satisfying the equation
\begin{eqnarray}\label{gd-eq}
0=\nabla_x f(\zb)u+\Big(\nabla\big(\nabla_x g(\cdot)^\ast\hat\lambda\big)(\zb)+{\cal H}_{\hat\lambda}\Big)(0,u)+\pi_{X^\ast}(w^\ast)=\F_{\hat\lambda}u +\pi_{X^\ast}(w^\ast)
\end{eqnarray}
with taking into account the notation in \eqref{H1} and \eqref{F}. Proceeding further as in the proof of Theorem~\ref{ThHoelder} with $w_k^\ast=w^\ast$, we arrive at a contradiction to \eqref{EqSuffCond}, which thus verifies \eqref{calm} in case (i).

Now we turn to the assumptions in (ii). Again arguing by contradiction as in the proof in case (i) gives us \eqref{sum} with some $u\ne 0$. Applying the graphical derivative computation in Theorem~\ref{CorStrictCompl} with $q=0\in\Pp$ therein, we find $z^\ast\in N_{T_{{\rm gph}\,\Gamma}(\zb)}(0,u)\cap\pi_{X^\ast}^{-1}(\xba)$, $\hat\lambda\in\Lambda(\zb,z^\ast)$, and $w^\ast\in N_{\K_{{\rm gph}\,\Gamma}(\zb,z^\ast)}(0,u)$ satisfying equation \eqref{gd-eq}. Proceeding then as in the proof of Theorem~\ref{ThHoelder} with $w_k^\ast=w^\ast$ leads us to $\skalp{\F_{\hat\lambda}u,u}=0$ and to $\skalp{\F_{\hat\lambda}u,\tilde u}\ge 0$ for all $\tilde u\in\K_{\Gamma(\pb)}$. Taking finally into account that $\hat\lambda\in\Lambda(\zb,\nabla g(\zb)^\ast\lambda)=\Lambda(\zb,\nabla g(\zb)^\ast\hat\lambda)$ and remembering the construction of $\Lb$ in \eqref{mult1} bring us to a contradiction with \eqref{EqStrongSuffCond} and thus complete the proof of the theorem.$\h$\vspace*{0.03in}

The next theorem shows that conditions \eqref{EqSuffCond} and \eqref{EqStrongSuffCond} are {\em necessary} for isolated calmness of PVS \eqref{EqSolMap} in fairly general settings.\vspace*{-0.07in}

\begin{Theorem}[\bf necessary conditions for isolated calmness in PVS]\label{calm-nec} Let PCS in \eqref{EqConstrSystem} enjoy the Robinson stability property at $(\pb,\xb)$, and let PVS in \eqref{EqSolMap} be isolatedly calm at this point. Impose further the following additional assumptions: $P=P_1\times P_2$,
\begin{eqnarray*}
g\big((p_1,p_2),x\big)=g\big((p_1,\pb_2),x\big)\;\mbox{ for all }\;\big((p_1,p_2),x\big)\in P\times X,\;\mbox{ and }\;\nabla_{p_2}f(\op,\ox)\;\mbox{ is surjective}.
\end{eqnarray*}
Then condition \eqref{EqSuffCond} is fulfilled being thus necessary for isolated calmness. If we have furthermore
\begin{eqnarray}\label{EqConstLambda1}
\begin{array}{ll}
\nabla h\big(g(\oz))\nabla_x^2g(\zb)(u,u)+\nabla^2h\big(g(\zb)\big)\big(\nabla_xg(\zb)u,\nabla_xg(\zb)u\big)\\
\in\rge\big(\nabla h(g(\zb))\circ\nabla g(\zb)\big)+\lin T_C\big(g(\zb)\big)\;\mbox{ for all}\;u\in\K_{\Gamma(\pb)}\big(\xb,-f(\zb)\big)
\end{array}
\end{eqnarray}
with $\oz=(\op,\ox)$, then condition \eqref{EqStrongSuffCond} is also satisfied.
\end{Theorem}\vspace*{-0.07in}
{\bf Proof.} To verify the first assertion, suppose that condition \eqref{EqSuffCond} fails and thus find $0\not=u\in\K_{\Gamma(\pb)}(\xb,\xba)$, $\lambda\in\Lb$, and $\hat\lambda\in\Lambda(\zb,\nabla g(\zb)^\ast\lambda;(0,u))$ such that $\skalp{\F_{\hat\lambda}u,u}\le 0$ and $\skalp{\F_{\hat\lambda}u,\tilde u}\ge 0$ for all $\tilde u\in\K_{\Gamma(\pb)}(\xb,-f(\zb))$. Taking $\tilde u=u$ yields $\skalp{\F_{\hat\lambda}u,u}=0$ and therefore
\begin{eqnarray*}
\big\la-\F_{\hat\lambda}u,\tilde u+\alpha u\big\ra\le 0\;\mbox{ whenever }\;\tilde u\in\K_{\Gamma(\pb)}(\xb,\xba),\alpha\in\R,
\end{eqnarray*}
which can be equivalently rewritten as
\begin{eqnarray*}
-\F_{\hat\lambda}u\in\big(\K_{\Gamma(\pb)}(\xb,\xba)+[u]\big)^\ast=\K_{\Gamma(\pb)}(\xb,\xba)^\ast\cap[u]^\perp.
\end{eqnarray*}
Thus there is an element $w_{x^\ast}\in\K_{\Gamma(\pb)}(\xb,\xba)^\ast\cap[u]^\perp$ satisfying
\begin{eqnarray*}
0=\F_{\hat\lambda}u+w_{x^\ast}=\nabla_x f(\zb)u+\Big(\nabla\big(\nabla_x g(\cdot)^\ast\hat\lambda\big)(\zb)+{\cal H}_{\hat\lambda}\Big)(0,u)+w_{x^\ast}.
\end{eqnarray*}
By Proposition~\ref{LemNormalCone} we get sequences $\lambda_k\in N_C(g(\xb))$ and $\alpha_k\in \R$ such that
\begin{eqnarray*}
\disp w_{x^\ast}=\lim_{k\to\infty}\nabla_xg(\zb)^\ast\lambda_k+\alpha_k\xba=\lim_{k\to\infty}\nabla_xg(\zb)^\ast(\lambda_k+\alpha_k\hat\lambda).
\end{eqnarray*}
Furthermore, defining $z^\ast:=\nabla g(\zb)^\ast\hat\lambda$ and $w_k^\ast:=\nabla g(\zb)^\ast(\lambda_k+\alpha_k\hat\lambda)=\nabla g(\zb)^\ast\lambda_k+\alpha_kz^\ast$ ensures by Proposition~\ref{LemNormalCone} again that $w_k^\ast\in\K_{{\rm gph}\,\Gamma}(\zb,z^\ast)$ together with
\begin{eqnarray*}
\disp\lim_{k\to\infty}\pi_{X^\ast}(w_k^\ast)=w_x^\ast,\;\lim_{k\to\infty}\skalp{w_k^\ast,(0,u)}=\lim_{k\to\infty}\skalp{\pi_{X^\ast}(w_k^\ast),u}
=\skalp{w_{x^\ast},u}=0,\;\mbox{ and }\;\hat\lambda\in\Lambda\big(\zb,z^\ast;(0,u)\big).
\end{eqnarray*}
Using now the graphical derivative estimate from Proposition~\ref{CorRemCases} with $q=0$ tells us that
\begin{eqnarray*}
0\in\nabla_x f(\zb)u+D\Psi(\zb,\xba)(0,u)=\nabla f(\zb)(0,u)+D\Psi(\zb,\xba)(0,u),
\end{eqnarray*}
and therefore there exist sequences $t_k\downarrow 0$, $u_k\to 0$, and $q_k\to 0$ such that
\begin{eqnarray*}
{\rm dist}\big(-f(\pb+t_kq_k,\xb+t_ku_k);\Psi(\pb+t_kq_k,\xb+t_ku_k)\big)=\oo(t_k)\;\mbox{ as }\;k\to\infty.
\end{eqnarray*}
It follows from the robustness of the imposed surjectivity/full rank condition on $\nabla_{p_2}f(\zb)$ that we can adjust the approximating sequence $\{q_k\}$ a bit so that
\begin{eqnarray*}
-f(\pb+t_k\tilde q_k,\xb+t_ku_k)\in\Psi(\pb+t_kq_k,\xb+t_ku_k)=\Psi(\pb+t_k\tilde q_k,\xb+t_ku_k)
\end{eqnarray*}
for $\tilde q_k:=q_k+(0,q_{2,k})$ with $q_{2,k}\st{P_2}{\to}0$ as $k\to\infty$. Hence $\xb+t_ku_k\in S(\pb+t_k\tilde q_k)$ and $0\not=u\in DS(\pb,\xb)(0)$, which tells us that $S$ is not isolatedly calm at $(\op,\ox)$, and thus the first assertion of the theorem is justified. To verify the necessary of condition \eqref{EqSuffCond} for the isolated calmness of $S$ in the second assertion, observe that the imposed additional requirement \eqref{EqConstLambda1} implies that
\begin{eqnarray*}
\lambda\in\Lambda(\zb,\nabla g(\zb)^\ast\lambda;(0,u))\;\mbox{ for all }\;u\in\K_{\Gamma(\pb)}\big(\xb,-f(\zb)\big)\;\mbox{ and all }\;\lambda\in\Lb,
\end{eqnarray*}
as can be checked similarly to the proof of Theorem~\ref{CorStrictComplAux}. Therefore conditions \eqref{EqSuffCond} and \eqref{EqStrongSuffCond} are equivalent in this case, and the proof of the theorem is complete. $\h$\vspace*{0.03in}

Finally,
we present an example demonstrating that isolated calmness of PVS \eqref{EqSolMap} may fail under conditions \eqref{EqSuffCond} and \eqref{EqStrongSuffCond}, respectively, when neither $\Pp=P$, nor $K$ is polyhedral, nor strict complementarity holds in second-order cone programming.\vspace*{-0.07in}

\begin{example}[\bf failure of isolated calmness of PVS generated by the ice-cream cone in $\R^3$]\label{calm-fail}
{\rm Consider the parametric variational system of type \eqref{EqSolMap} given by
\begin{eqnarray*}0&\in&\Big(-\sqrt{2}+\frac{\sqrt{2}}4(2x_1+p_1+p_2)-\xi(p_1-p_2), 1+\frac{\sqrt{2}}4x_2\Big)+\widehat N_{\Gamma(p_1,p_2)}(x_1,x_2),\\
\Gamma(p_1,p_2)&:=&\big\{(x_1,x_2)\big|\;g(p,x):=(x_1+p_1,x_1+p_2,x_2)\in{\cal Q}_3\big\},\\
\xi(t)&:=&\sqrt{2}\Big(\big(1+\vert t\vert^{\frac 23}\big)^{-\frac 12}+\frac 12\vert t\vert^{\frac 23}\Big),
\end{eqnarray*}
where ${\cal Q}_3:=\{y\in\R^3\mv y_3\ge\sqrt{y_1^2+y_2^2}\}$ denotes the ice-cream cone in $\R^3$. It is easy to check that the function $\xi$ is continuously differentiable with $\xi'(0)=0$ and that $S(p)$ consists of the points satisfying the basic first-order optimality conditions for the parametric {\em quadratic second-order cone program}:
\begin{eqnarray*}
&&\min_{x_1,x_2}\frac{\sqrt{2}}4 x_1^2+\frac{\sqrt{2}}8 x_2^2+\Big(\frac{\sqrt{2}}4(p_1+p_2)-\sqrt{2}-\xi(p_1-p_2)\Big)x_1 +x_2\\
&&\mbox{subject to }\;(x_1+p_1,x_1+p_2,x_2)\in{\cal Q}_3.
\end{eqnarray*}
The cost function of this program is obviously strictly convex for every $p$, and the classical Slater constraint qualification holds at every feasible point. Thus for every $p\in\R^2$ the solution set $S(p)$ is a singleton. In fact, the function $\xi$ was chosen in such a way that
\begin{eqnarray*}
S(p)=\left\{\Big(\frac 12\big(-p_1-p_2+\vert p_1-p_2\vert^{\frac 23}),\frac{\vert p_1-p_2\vert^{\frac 23}}{\sqrt{2}}\sqrt{1+\vert p_1-p_2\vert^{\frac 23}}\Big)\right\}.
\end{eqnarray*}
Considering the reference pair $(\op,\ox)$ with $\pb=(0,0)$ and $\xb=(0,0)$, we see that $S$ is not isolatedly calm at $(\pb,\xb)$. Observe that the standing assumptions are satisfied with $h(y)=y$ and $K={\cal Q}_3$, which is not polyhedral. By the Slater constraint qualification, the Robinson stability condition is also fulfilled for the system $g(p,x)\in C$ at $(\pb,\xb)$. The parameter characteristic subspace here is $\Pp=\{(\alpha,\alpha)\mv\alpha\in\R\}$, and so  $\Pp\not=P$. Furthermore, the multipliers set \eqref{mult1} is $\Lb=\big\{(\frac 1{\sqrt{2}},\frac 1{\sqrt{2}},-1)\big\}$, which shows that the strict complementarity condition fails at the reference point as well.

On the other hand, we have $\F_\lambda u=\big(\frac{\sqrt{2}}2u_1,\frac{\sqrt{2}}4u_2\big)$, and thus condition \eqref{EqSuffCond} is satisfied. Applying Theorem~\ref{ThHoelder} gives us the H\"olderian isolated calmness
\begin{eqnarray*}
S(p)\cap U\subset\big\{\xb\big\}+\ell\big(\norm{p-\pb}+\sqrt{{\rm dist}(p-\pb;\Pp)}\big)\B_{\R^2}\;\mbox{ for all }\;p\in Q
\end{eqnarray*}
with some modulus $\ell>0$ and neighborhoods $U$ of $\xb$ and $Q$ of $\pb$. In fact, the distance expression ${\rm dist}(p-\pb;\Pp)=\frac 1{\sqrt{2}}\vert p_1-p_2\vert$ brings us to an upper estimate of the entire solution set $S(p)$ given by
\begin{eqnarray*}
S(p)\subset\big\{\xb\}+\ell(\norm{p-\pb}+{\rm dist}(p-\pb;\Pp)^{\frac 23})\B_{\R^2},\quad p\in Q.
\end{eqnarray*}}
\end{example}\vspace*{-0.2in}

\section{Concluding Remarks}\vspace*{-0.1in}

This paper demonstrates that applying advanced tools and techniques of second-order variational analysis allows us to efficiently compute the graphical derivative of the normal cone mapping generated by a large class of parametric constraint systems while providing in this way striking new developments even for nonparametric nonpolyhedral systems without any nondegeneracy and/or metric regularity assumptions. The established results are applied to deriving necessary and sufficient conditions for isolated calmness and its H\"olderian counterpart in parametric variational systems that arise, in particular, as first-order optimality conditions in parametric conic programming.

The obtained second-order calculations have strong potentials for further applications to optimization and related areas where the usage of graphical derivatives for normal cone mappings is highly beneficial. To this end, we mention the recent developments in \cite{MS17} on applying the graphical derivative computation to the study of {\em critical multipliers} in polyhedral variational systems that are largely responsible for slow convergence of primal-dual algorithms in optimization. It seems that employing the new results obtained here would lead us to a significant progress in this direction.\vspace*{0.05in}

{\bf Acknowledgements.} The research of the first author was partially supported by the Austrian Science Fund (FWF) under grant P29190-N32. The research of the second author was partially supported by the USA National Science Foundation under grant DMS-1512846, by the USA Air Force Office of Scientific Research under grant No.\,15RT0462, and by the RUDN University Program 5-100.
\end{document}